\documentclass[reqno,11pt]{amsart}

\usepackage[a4paper,left=23mm,right=23mm,top=30mm,bottom=30mm,marginpar=25mm]{geometry}
\usepackage{amsmath}
\usepackage{amssymb}
\usepackage{amsthm}
\usepackage{amscd}
\usepackage{mathtools}
\usepackage{url}
\usepackage[ansinew]{inputenc}
\usepackage{color}
\usepackage[final]{graphicx}
\usepackage{enumitem}
\usepackage{esint}
\usepackage{bbm}
\usepackage{subfigure}

\usepackage{tikz}
\usetikzlibrary{arrows,shapes,automata,backgrounds,petri,positioning}
\usetikzlibrary{decorations.pathmorphing}
\usetikzlibrary{decorations.shapes}
\usetikzlibrary{decorations.text}
\usetikzlibrary{decorations.fractals}
\usetikzlibrary{decorations.footprints}
\usetikzlibrary{shadows}
\usetikzlibrary{calc}
\usetikzlibrary{spy}

\usepackage{cite}
\usepackage[hidelinks]{hyperref}

\renewcommand{\epsilon}{\varepsilon}

\numberwithin{equation}{section}

\newtheoremstyle{thmlemcorr}{10pt}{10pt}{\itshape}{}{\bfseries}{.}{10pt}{{\thmname{#1}\thmnumber{ #2}\thmnote{ (#3)}}}
\newtheoremstyle{thmlemcorr*}{10pt}{10pt}{\itshape}{}{\bfseries}{.}\newline{{\thmname{#1}\thmnumber{ #2}\thmnote{ (#3)}}}
\newtheoremstyle{defi}{10pt}{10pt}{\itshape}{}{\bfseries}{.}{10pt}{{\thmname{#1}\thmnumber{ #2}\thmnote{ (#3)}}}
\newtheoremstyle{remexample}{10pt}{10pt}{}{}{\bfseries}{.}{10pt}{{\thmname{#1}\thmnumber{ #2}\thmnote{ (#3)}}}
\newtheoremstyle{ass}{10pt}{10pt}{}{}{\bfseries}{.}{10pt}{{\thmname{#1}\thmnumber{ A#2}\thmnote{ (#3)}}}

\theoremstyle{thmlemcorr}
\newtheorem{theorem}{Theorem}
\numberwithin{theorem}{section}
\newtheorem{lemma}[theorem]{Lemma}
\newtheorem{corollary}[theorem]{Corollary}
\newtheorem{proposition}[theorem]{Proposition}

\theoremstyle{thmlemcorr*}
\newtheorem{theorem*}{Theorem}
\newtheorem{lemma*}[theorem]{Lemma}
\newtheorem{corollary*}[theorem]{Corollary}
\newtheorem{proposition*}[theorem]{Proposition}
\newtheorem{problem*}[theorem]{Problem}
\newtheorem{conjecture*}[theorem]{Conjecture}

\theoremstyle{defi}
\newtheorem{definition}[theorem]{Definition}

\theoremstyle{remexample}
\newtheorem{remarkx}[theorem]{Remark}
\newenvironment{remark}
  {\pushQED{\qed}\remarkx}
  {\popQED\endremarkx}
\newtheorem{examplex}[theorem]{Example}
\newenvironment{example}
  {\pushQED{\qed}\examplex}
  {\popQED\endexamplex}

\theoremstyle{ass}

\newcommand{\Acal}{\mathcal{A}}

\newcommand{\Ccal}{\mathcal{C}}
\newcommand{\Dcal}{\mathcal{D}}

\newcommand{\Fcal}{\mathcal{F}}

\newcommand{\Hcal}{\mathcal{H}}

\newcommand{\Lcal}{\mathcal{L}}

\newcommand{\Pcal}{\mathcal{P}}

\newcommand{\Qcal}{\mathcal{Q}}

\newcommand{\Rcal}{\mathcal{R}}
\newcommand{\Scal}{\mathcal{S}}

\DeclareMathOperator{\dist}{dist}

\DeclareMathOperator{\supp}{supp}

\newcommand{\abslr}[1]{\left|#1\right|}
\newcommand{\abss}[1]{\left\langle #1 \right\rangle}

\newcommand{\absb}[1]{\bigl|#1\bigr|}

\newcommand{\N}{\mathbb{N}}
\newcommand{\R}{\mathbb{R}}

\newcommand{\C}{\mathbb{C}}

\newcommand{\weakly}{\rightharpoonup}

\newcommand{\eps}{\epsilon}

\DeclareMathOperator{\Div}{div}

\newcommand{\D}{\mathcal{D}}

\def\XXint#1#2#3{{\setbox0=\hbox{$#1{#2#3}{\int}$}
\vcenter{\hbox{$#2#3$}}\kern-.5\wd0}}

\DeclareMathOperator{\Id}{Id}

\definecolor{KUblue}{RGB}{0,102,153}

\DeclarePairedDelimiter\abs{\lvert}{\rvert}
\DeclarePairedDelimiter{\norm}{\lVert}{\rVert}
\DeclarePairedDelimiter{\inner}{\langle}{\rangle}

\newcommand{\Rn}{\R^{n}}

\renewcommand{\phi}{\varphi}

\newcommand{\starto}{\stackrel{*}{\rightharpoonup}}
\newcommand{\weakto}{\rightharpoonup}

\newcommand{\Trho}{T_{\rho(\cdot)}}

\newcommand{\Srho}{S_{\rho(\cdot)}}

\def\XXint#1#2#3{{\setbox0=\hbox{$#1{#2#3}{\int}$}
     \vcenter{\hbox{$#2#3$}}\kern-.5\wd0}}

\newcommand{\Rmn}{\mathbb{R}^{N \times n}}

\makeatletter
\g@addto@macro\bfseries{\boldmath}
\makeatother

\usepackage{esint}

\usepackage{appendix}

\renewcommand{\O}{\Omega}

\renewcommand{\d}{\delta}

\newcommand{\Gr}{D_{\rho(\cdot)}}
\newcommand{\Hr}{H^{\rho(\cdot),p}(\Omega)}
\newcommand{\Hrm}{H^{\rho(\cdot),p}(\Omega;\R^N)}
\newcommand{\Hrg}{H^{\rho(\cdot),p}_g(\Omega)}
\newcommand{\Hro}{H^{\rho(\cdot),p}_0(\Omega)}
\newcommand{\Hrmg}{H^{\rho(\cdot),p}_g(\Omega;\R^N)}
\newcommand{\Hrr}{H^{\rho(\cdot),p}(\R^n)}
\newcommand{\Nr}{N^{\rho(\cdot),p}(\Omega)}
\newcommand{\Smm}{S^m_{\nu,\mu}}

\newcommand{\Fr}{\Fcal_{\rho(\cdot)}}
\newcommand{\Op}{\mathrm{Op}}
\newcommand{\OPS}{\mathrm{OPS}}
\newcommand{\OPSmm}{\OPS^{m}_{\nu,\mu}}
\newcommand{\LBH}{\mathrm{LBH}}

\newcommand{\Qr}{\Qcal_{\rho(\cdot)}}
\newcommand{\Qro}{\Qcal_{\rho(\cdot)}^{\Omega}}
\newcommand{\qr}{q_{\rho(\cdot)}}
\newcommand{\Cr}{\Ccal_{\rho(\cdot)}}
\newcommand{\Cro}{\Ccal_{\rho(\cdot)}^{\Omega}}
\newcommand{\PR}{\Pcal_{\rho(\cdot)}}
\newcommand{\PROmega}{\Pcal_{\rho(\cdot), \Omega}}
\newcommand{\pr}{p_{\rho(\cdot)}}
\newcommand{\Rr}{\Rcal_{\rho(\cdot)}}
\newcommand{\RrOmega}{\Rcal_{\rho(\cdot), \Omega}}
\newcommand{\Rrr}{\Rcal'_{\rho(\cdot)}}
\newcommand{\RrrOmega}{\Rcal'_{\rho(\cdot), \Omega}}
\newcommand{\Grr}{\overline{D}_{\rho(\cdot)}}
\newcommand{\Er}{E_{\rho(\cdot)}}

\newcommand{\black}{\color{black}}

\newcommand{\Cinfo}{C_0^{\infty}(\overline{\Omega})}
\newcommand{\Cinf}{C^{\infty}(\overline{\Omega})}

\usepackage{amsmath}
\usepackage{latexsym}
\usepackage{amssymb}

\def\XXint#1#2#3{{\setbox0=\hbox{$#1{#2#3}{\int}$}
		\vcenter{\hbox{$#2#3$}}\kern-.5\wd0}}

\usepackage{graphicx}
\usepackage{wrapfig}

\title[Local boundary conditions in nonlocal hyperelasticity]{Local boundary conditions in nonlocal hyperelasticity via heterogeneous horizons}

\author{Carolin Kreisbeck}
\address{Mathematisch-Geographische Fakult\"at, Katholische Universit\"at Eichst\"att-Ingolstadt, Osten\-stra\-{\ss}e 28, 85072 Eichst\"att, Germany}
\email{carolin.kreisbeck@ku.de}

\author{Hidde Sch\"{o}nberger}
\address{Research Institute in Mathematics and Physics, UCLouvain, Chemin du Cyclotron 2, 1348 Louvain-la-Neuve, Belgium}
\email{hidde.schonberger@uclouvain.be}

\begin{document}
\maketitle

\thispagestyle{empty}
\begin{abstract}  

 In this paper, we consider a class of variational problems with integral functionals involving nonlocal gradients. These models have been recently proposed as refinements of classical hyperelasticity, aiming for an effective framework to capture also discontinuous and singular material effects. 
Specific to our set-up is a space-dependent interaction range that vanishes at the boundary of the reference domain. This ensures that the nonlocal operator depends only on values within the domain and localizes to the classical gradient at the boundary, which allows for a seamless integration of nonlocal modeling with local boundary values. The main contribution of this work is a comprehensive theory for the newly introduced associated Sobolev spaces, including the rigorous treatment of a trace operator and Poincar\'e inequalities. A central aspect of our technical approach lies in exploiting connections with pseudo-differential operator theory.  As an application, we establish the existence of minimizers for functionals with quasiconvex or polyconvex integrands depending on heterogeneous nonlocal gradients, subject to local Dirichlet, Neumann or mixed-type boundary conditions.

\vspace{8pt}

%
\vspace{8pt}

\vspace{8pt}

\noindent\textsc{Date:} \today.
\end{abstract}


\maketitle
\thispagestyle{empty}

\section{Introduction}
Nonlocal modeling, where the state at a point depends not only on its immediate neighborhood but also on an interacting region of finite horizon or even the entire domain, extends beyond classical differential formulations. This approach naturally gives rise to equations or variational problems involving integral, integro-differential, fractional, or other nonlocal operators. Notable benefits over classical local models show especially in situations that require higher accuracy, greater generality, or increased robustness. 
Within nonlocal frameworks, one can account for global dependencies and long-range interactions, and thereby bridge scales, while their weaker regularity requirements make them particularly well-suited for capturing discontinuities, singularities and other complex phenomena. 
As a result, these models remain valid in regimes where classical ones tend to fail, such as in fracture, plasticity, and heterogeneous media, to mention just a few examples from solid mechanics. 
 However, these advantages come at a cost: Nonlocal models are typically computationally demanding, and a sound treatment of boundary conditions and interfaces, which are inherently local, is a
challenge. In both regards, classical local models remain clearly superior.

With the aim of combining the advantages of both modeling paradigms in a mathematically tractable framework, local-to-nonlocal coupling seeks to achieve the best of both worlds. It integrates the higher accuracy, reduced regularity requirements, and nonlocal features available within the domain with the simplicity of well-defined local boundary or interface conditions. 
The literature on local-to-nonlocal coupling mechanisms (see for instance the review article \cite{DLSTY19}, as well as the references therein) comprises a variety of strategies, which can be classified basically into two categories: constant-horizon and varying-horizon approaches. Roughly speaking, the first class of methods introduces abrupt changes in the horizon, relying on domain decomposition in which different subdomains described by local and nonlocal models are coupled through transmission conditions or overlapping regions. In contrast, the second class -- which is the focus of this work -- shrinks the interaction range near boundaries or interfaces to obtain a smooth transition to a local description.
Even though this general idea of coupling through heterogeneous localization may seem intuitive, establishing it rigorously for a specific set-up requires formulating a tailored functional analytic framework and developing the corresponding technical tools. 
In the setting of fractional Sobolev spaces defined via Gagliardo-type seminorms and with motivation from peridynamics~\cite{Sil00, SLS14}, such a study was initiated in \cite{TiD17} and further developed in \cite{TTD19, DMT22, DTWY22} as well as in the recent papers \cite{ScD24, ScD25}.

In this work, we address coupling via heterogeneous localization in the context of the recently introduced nonlocal hyperelasticity \cite{BeCuMC, BeCuMC23}, which can be interpreted as falling under the theory of state-based peridynamics \cite{SEWXA07}. In comparison with standard hyperelastic models, the classical deformation gradient is replaced by a nonlocal analogue -- an averaged linear approximation that accounts for interactions within the horizon. The key objects in these models are thus nonlocal gradients, operators that have attracted considerable attention in the analysis community in recent years. The novel aspect we study here is the effect of a spatially varying horizon, with the aim of better understanding and establishing a theory for heterogeneous nonlocal gradients and their associated Sobolev function spaces (see Section~\ref{subsec:intro_heterogeneousgradients}). One of our central goals is to rigorously establish a local trace theorem, showing that functions in our heterogeneous nonlocal Sobolev spaces are generally less regular inside the domain, but still admit well-defined local boundary values.
Let us note that, while the focus in the following lies on boundary conditions, the same technical approach is applicable to interfaces as well. 

\subsection{From homogeneous to heterogeneous nonlocal gradients}\label{subsec:intro_heterogeneousgradients}
At the core of this work
 are space-dependent variants of nonlocal gradients, which we introduce here as follows:
 Given a bounded and connected Lipschitz domain $\Omega\subset \mathbb{R}^n$, a smooth 
 horizon function $\delta:\Omega\to (0, \infty)$ with $\delta(x)\leq \dist(x, \partial \Omega)$, and a radial kernel function $\rho:\R^n\setminus \{0\}\to [0, \infty)$, the heterogeneous nonlocal gradient $\Gr$ is the singular integral operator defined for a function $\varphi\in C^\infty(\overline{\Omega})$ by
	\begin{align}\label{def:Drhoi_intro}
		D_{\rho(\cdot)} \varphi (x) =\delta(x)^{-n} \int_{\Omega} \frac{\varphi (y)- \varphi (x)}{|y-x|}\frac{y-x}{|y-x|} \rho \Bigl(\frac{y-x}{\delta(x)}\Bigr)\,dy \quad\text{ for $x \in \Omega$.}
	\end{align} 
	
To be more precise about the assumptions, we suppose that the function $\delta$ is smooth and decays
  sufficiently fast toward the boundary. Moreover, the kernel $\rho$, 
 which is assumed to be supported on the closed unit ball of $\mathbb{R}^n$, i.e.,~$\supp \rho \subset B_1(0)$, falls within the general framework of kernels for nonlocal gradients with constant horizon introduced in~\cite{BeMCSc24} (see 
Section~\ref{subsec:homogeneous_nonlocalgradients} below
 for 
 more context and examples). In particular, this imposes in addition to technical conditions regarding smoothness, integrability and monotonicity, the hypothesis that 
$\rho$ lies, near the origin, in between two fractional kernels of orders $\lambda\in (0,1)$ and $\kappa\in (0, 1)$ with $\lambda\leq \kappa$; for a precise formulation of these assumptions, we refer to~\ref{itm:h0}-\ref{itm:upper} in Section~\ref{sec:nonlocalgradients}.

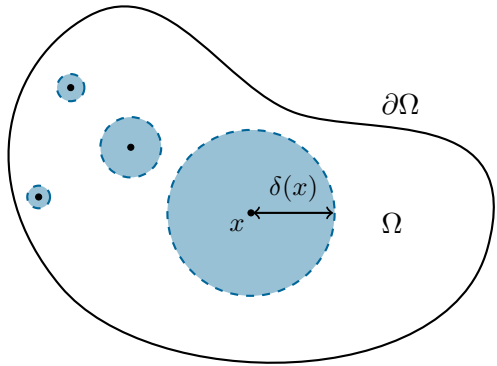
\begin{figure}[h!]
\centering

\begin{tikzpicture}[x=0.75pt,y=0.75pt,yscale=-1,xscale=1, spy using outlines]

\draw [thick] (244,56) .. controls (284,34) and (312,93) .. (348,105) .. controls (384,117) and (463,100) .. (443,173) .. controls (423,246) and (275,246) .. (230,192) .. controls (185,138) and (204,78) .. (244,56) -- cycle ;

\draw (385,150) node [anchor=north west] {$\Omega$};
\draw (385,90) node [anchor=north west] {$\partial\Omega$};
\draw [color = KUblue, thick, dashed,  fill = KUblue, fill opacity =0.4] (325,155) circle [radius = 1.1 cm];
\draw [color = black, fill=black] (325,155) circle [radius = 0.04cm];
\draw (327,153) node [anchor = north east] {{\small $x$}};
\draw [color = KUblue, thick, dashed,  fill = KUblue, fill opacity =0.4] (265,122) circle [radius = 0.4cm];
\draw [color = black, fill=black] (265,122) circle [radius = 0.04cm];
\draw [color = KUblue, thick, dashed,  fill = KUblue, fill opacity =0.4] (235,92) circle [radius = 0.18cm];
\draw [color = black, fill=black] (235,92) circle [radius = 0.04cm];

\draw [color = KUblue, thick, dashed,  fill = KUblue, fill opacity =0.4] (219,147) circle [radius = 0.15cm];
\draw [color = black, fill=black] (219,147) circle [radius = 0.04cm];
\draw[<->,thick] (326,155) -- (366,155);
\draw (329,131) node [anchor=north west] {{\small $\delta(x)$}};

\end{tikzpicture}
\caption{Illustration of the reference configuration $\Omega$ with the space-dependent horizon $\delta(x)$ for selected points $x$ and the local boundary $\partial \Omega$.}
\label{fig:heterogeneous_horizon}
\end{figure}

The definition of $\Gr$ in~\eqref{def:Drhoi_intro} reflects that the heterogeneous nonlocal gradient arises from a homogeneous nonlocal gradient via  a position-dependent rescaling involving the horizon $\delta(x)$ for $x\in \Omega$. Accordingly, the nonlocal interaction range for computing $\Gr$ of a function at a point $x\in \Omega$ is given by $\delta(x)$, cf.~Figure~\ref{fig:heterogeneous_horizon}. 
 The motivation behind this intuitive approach is to recover the classical local gradient at the boundary $\partial \Omega$, due to the vanishing of the horizon in the limit $\dist(x,\partial\Omega)\to 0$.

After suitably extending the notion of heterogeneous nonlocal gradient from smooth functions to larger classes, such as locally integrable functions, and more generally, distributions,
which is a central aspect in this work (see Section~\ref{subsec:defnonlocalgradient} for more details), 
one can introduce heterogeneous nonlocal Sobolev spaces. For $p\in (1, \infty)$, these are defined in analogy to the standard Sobolev spaces as
\begin{align*}
\Hr = \{u\in L^p(\Omega): \Gr u\in L^p(\Omega;\R^n)\},
\end{align*}
equipped with the norm $\norm{u}_{\Hr} =\norm{u}_{L^p(\Omega)} + \norm{\Gr u}_{L^p(\Omega;\R^n)}$.

\subsection{Background on nonlocal gradients with constant horizon}\label{subsec:homogeneous_nonlocalgradients} 
 As indicated above, our results in the setting with spatially varying horizons build on the theory of homogeneous nonlocal gradients   as a foundation. 
 Let us briefly review 
 some developments concerning these operators, namely, $D_\rho$ with a kernel function $\rho$ defined for $\varphi \in C_c^\infty(\R^n)$ by 
	\begin{align*}
		D_{\rho} \varphi (x) =\ \int_{\R^n} \frac{\varphi (y)- \varphi (x)}{|y-x|}\frac{y-x}{|y-x|} \rho (y-x)\,dy \quad\text{ for $x \in \R^n$.}
	\end{align*}  
 
 A prototypical example of a nonlocal gradient that has been intensively studied in recent years is the Riesz fractional gradient $D^s$ for fractional order $s\in (0,1)$, whose kernel is a multiple of the Riesz potential of order $1-s$, see \cite{ShS2015, ShS2018}. 
 It represents a natural choice of fractional derivative, as it combines desirable mathematical properties, such as its distributional character~\cite{Comi1, Comi2}, with structural features like homogeneity and translation invariance \cite{Silhavy2019}, which are particularly relevant for modeling purposes.

Motivated by applications in continuum mechanics, which necessitate working on bounded domains,   \cite{BeCuMC23} introduced and studied a finite-horizon version $D_\delta^s$ of the Riesz fractional gradient $D^s$ with horizon $\delta>0$, obtained by truncation of the potential; some basic properties and technical tools for $D_\delta^s$ were further investigated also in~\cite{CuKrSc23}.
 Extending and unifying these special cases, the authors of~\cite{BeMCSc24} established a generalized framework of nonlocal gradients $D_\rho$ with suitable properties on the kernel $\rho$, which we use as a basis in this paper.
For the function spaces associated to the general nonlocal gradients $D_\rho$,
 denoted by $H^{\rho, p}(\R^n)$ for $p\in (1, \infty)$, a toolbox analogous to that of the classical Sobolev space theory has been established in the references mentioned above. Important ingredients include continuous and compact embeddings, theorems on the density of smooth functions, and estimates such as the Poincar\'e inequality. 
A notable observation about the connection between nonlocal and classical gradients gives rise to a translation mechanism \cite{CuKrSc23, CKS25, KrS22}, which has turned out to be a useful method for transferring results from the local to the nonlocal setting. Essentially, it states the existence of a bounded, linear convolution operator $\Qcal_\rho: H^{\rho, p}(\R^n)\to W^{1,p}(\R^n)$ with integrable potential and inverse $\Pcal_\rho$ such that 
\begin{align}\label{intro:translation_class}
\qquad D_\rho = \nabla \Qcal_\rho \qquad \text{and} \qquad \nabla = D_\rho \Pcal_\rho. 
\end{align} 
Note that in the classical fractional case $D_
\rho=D^s$, these translation operators correspond to the convolution with the Riesz potential of order $1-s$ and the fractional Laplacian $(-\Delta)^{(1-s)/2}$, respectively, cf.~also~\cite{ShS2015}. 
Another powerful tool to point out is the nonlocal counterpart of a fundamental theorem of calculus in \cite{BeCuMC23, BeMCSc24}, which, in short, allows one to recover a function from its nonlocal gradient via convolution with a suitable kernel function. 

Among recent results on variational problems in the context of nonlocal hyperelasticity
one can mention, for instance, the study of weak lower semicontinuity of nonlocal energy functionals of integral type \cite{KrS22, CuKrSc23}, various localization theorems (including in the limits of vanishing horizon $\delta\to 0$ and fractional parameter tending to 
$1$) \cite{MeS, CuKrSc23, CKS25, BCM21}, a characterization of zero nonlocal gradients and their relation to nonlocal Neumann-type problems \cite{KrS24, BCFR24}, general $\Gamma$-convergence results implying, among others, relaxation and homogenization formulas \cite{CuKrSc23, CCM25}, as well as a novel fractional modeling approach for dislocations in strain-gradient plasticity \cite{ACFS25}.

Beyond the line of research discussed above, one can find a range of works on nonlocal gradients; see, e.g., \cite{DGLZ, Delia, HaT23, HMT24, Arr25, CaR23} for alternative formulations and perspectives in nonlocal vector calculus.

\medskip

\color{black}
\subsection{Main results: Nonlocal Sobolev spaces}\label{subsec:intro_main1}  
After introducing nonlocal gradients with space-dependent horizons $\Gr$ as new operators, the task is to study their fundamental properties and develop a comprehensive theory for the associated heterogeneous nonlocal Sobolev spaces $\Hr$, which constitutes the main contribution of this work. 
We establish a technical toolbox, analogous to that available for classical Sobolev spaces, as a foundation for addressing problems in the context of partial differential equations and the calculus of variations where these spaces provide the natural class of admissible functions. 
In summary, here is an overview of our results; for full details, see Section~4:
\smallskip

\begin{itemize}
\color{black}
\item[\textit{(a)}] \textit{Continuous and compact embeddings:} To relate 
 our new heterogeneous nonlocal Sobolev spaces to classical spaces from the literature, we show that $\Hr$ lies continuously embedded in between the classical Sobolev space and the Bessel potential space $H^{\lambda, p}(\Omega)$, with the latter compactly contained in $L^p(\Omega)$; in formulas, 
 \begin{align*}
W^{1,p}(\Omega) \hookrightarrow \Hr \hookrightarrow H^{\lambda,p}(\Omega) \hookrightarrow\hookrightarrow L^p(\Omega).
\end{align*}  

\item[\textit{(b)}] \textit{A translation mechanism up to lower-order operators:} 
While --  in contrast to the constant-horizon setting --  a perfect translation is not available here, we can still derive a rigorous connection between heterogeneous nonlocal gradients and classical local ones that serves as a powerful  tool in our analysis. 
\color{black} In fact, we establish that there are bounded, linear operators 
\begin{center}
$\Qr^\Omega:\Hr \to W^{1,p}(\Omega)$ \quad and  \quad $\PROmega:W^{1,p}(\Omega) \to \Hr$ 
\end{center}
with the properties that 
\begin{align*}
\quad \ \ \qquad \Gr \approx \nabla \Qr^\Omega \qquad \text{and} \qquad \nabla \approx \Gr \PROmega,
\end{align*}
where $\approx$ refers to equality up to operators of lower order, which are manageable in the  applications of our interest; for a precise formulation of the statement, see Lemma~\ref{le:connection}. 
As a direct consequence, we also deduce the following intuitive 
 characterization of the space $\Hr$ as
\begin{align*}
\Hr=\{u\in L^p(\Omega): \Qr^\Omega u\in W^{1,p}(\Omega)\},
\end{align*}
with the norm $\norm{\cdot}_{\Hr}$ equivalent to $\norm{\Qr^\Omega\,\cdot\,}_{W^{1,p}(\Omega)}$.  
\smallskip
\item[\textit{(c)}]  \textit{Existence of a trace operator and classical boundary values:}
Despite its nonlocal nature, the space $\Hr$ has the characteristic feature that its functions possess a trace, and even share the same boundary values as classical Sobolev functions in $W^{1,p}(\Omega)$. This property is rooted in this statement:
There exists a unique bounded linear operator $\Trho:\Hr \to W^{1-1/p,p}(\partial \Omega)$ such that
\begin{align*}
\Trho u = T_{W^{1,p}(\Omega)} u \qquad \text{for all $u \in W^{1,p}(\Omega)$},
\end{align*}
and $\Trho$ possesses a continuous right-inverse; here, $T_{W^{1,p}(\Omega)}$ stands for the standard trace operator on the Sobolev space $W^{1,p}(\Omega)$.
Based on this trace theorem, one can introduce nonlocal spaces with prescribed classical boundary values, which we denote by $H_g^{\rho(\cdot), p}(\Omega)$ for $g\in W^{1-1/p, p}(\partial \Omega)$. \smallskip

\color{black}
\item[\textit{(d)}]  \textit{Density of smooth functions:} As a nonlocal analog of the Meyers-Serrin theorem for classical Sobolev spaces, we establish that the spaces of smooth functions $\Cinf$ and  $C_c^{\infty}(\Omega)$ are dense in $\Hr$ and $\Hro$, respectively, or expressed equivalently,
\begin{align*}
\Hr = \overline{\Cinf}^{\Hr} \qquad \text{and} \qquad \Hro = \overline{C_c^\infty(\Omega)}^{\Hr}. 
\end{align*}

\item[\textit{(e)}] \textit{Extension to functions on $\R^n$:}
We construct a linear and bounded operator $\Er$ that is well-suited for extending functions in $\Hr$ to elements in its corresponding full-space version $\Hrr$, see Definition~\ref{def:newspacesRn}. The approach is based on the Rychkov universal extension operator, combined with the translation operator $\Qr^\Omega$. 
 \smallskip

\item[\textit{(f)}]  \textit{Regularity results:} Our analysis confirms that the spaces $\Hr$ also contain functions less regular than Sobolev functions. In particular, when $\kappa p < n$ and $\kappa p < 1$, one can observe cavitation effects and jumps, respectively. Indeed, we obtain that for any open set $U \subset \R^n$ compactly contained in $\Omega$,
\[
H^{\kappa,p}_0(U)|_{\Omega} + W^{1,p}(\O) \subset \Hr. 
\]
Note that in the special case with fractional orders
 $\kappa=s=\lambda$, this inclusion along with the embedding in \textit{(a)} implies an identification of $\Hr$ with $H^{s,p}(\Omega)$ away from the boundary $\partial \Omega$.
\smallskip

\item[\textit{(g)}]  \textit{Characterization of functions with zero heterogeneous nonlocal gradient:} 
Identifying functions with zero gradient is far from trivial in the nonlocal setting. This is demonstrated, for example, in \cite{KrS24} by showing that the homogeneous finite-horizon fractional gradients may have an infinite-dimensional kernel.  We show here, in contrast, that the kernel  of $\Gr$ is always finite-dimensional.  Moreover, this kernel consists only of constant functions, and thus, mirrors the behavior of the classical local gradient, provided the horizon function takes sufficiently small values. 
Whether the same statement holds true without this additional technical assumption required in our approach remains open at this point.

\smallskip
\item[\textit{(h)}]  \textit{Poincar{\'e} inequalities in different versions:} 
Given that Poincar\'e inequalities are a fundamental tool of Sobolev space theory and play a crucial role in deriving compactness results, they are essential in our setting as well. For the heterogeneous nonlocal Sobolev spaces, 
we derive both a Poincar\'e-Wirtinger-type inequality and a variant for functions with partial zero boundary conditions,  building on the insights from \textit{(g)}. Specifically, we prove that there exists a constant $C>0$ such that
\[
\norm{u}_{L^p(\Omega)} \leq C \norm{\Gr u}_{L^p(\Omega;\R^n)}
\]
for all $u\in \Hr$ satisfying either $\int_\Omega u\,dx =0$ or $\Trho u=0$ $\Hcal^{n-1}$-a.e.~on some Borel-measurable set $\Gamma \subset \partial \Omega$ with $\Hcal^{n-1}(\Gamma) >0$.
\end{itemize}

The proofs of these results, and more broadly our underlying technical approach, rely crucially on a connection with classical pseudo-differential operator theory (see~e.g.~\cite{Kum81,Hor07}), which we establish at the outset and systematically exploit throughout.
More precisely, our analysis builds on the fundamental observation that the heterogeneous nonlocal gradient  $D_{\rho{(\cdot)}}$ corresponds to a (restricted) pseudo-differential operator with H\"ormander symbol. 
The existence of dual operators then allows us to introduce a distributional definition of  $D_{\rho{(\cdot)}}$ and to derive an integration-by-parts formula. We also employ the mapping properties of pseudo-differential operators between Bessel potential spaces, and the construction of parametrices as almost inverses plays an important role in defining the translation operator $\PROmega$ in \textit{(b)};  for comparison, in the homogeneous counterpart in~\eqref{intro:translation_class} one works in Fourier space to obtain $
\Pcal_\rho$ by inverting the convolution operator $\Qcal_\rho$.
As a final example, let us mention \textit{(g)}, where we rely on insights into the (partial) invertibility of  pseudo-differential operators with slowly varying symbols, which apply to $D_{\rho{(\cdot)}}$ when the horizon is sufficiently small.

Overall, considering that pseudo-differential operator theory  extends Fourier analysis by allowing for space-dependent symbols, the techniques in this paper provide a natural generalization of the reasoning used for homogeneous nonlocal gradients, e.g.~in~\cite{BeMCSc24, CuKrSc23, BeCuMC23}, to the heterogeneous setting.

\subsection{Main results: Applications to nonlocal variational problems.}\label{subsec:intro_main2}  With the above theoretical framework in place, we are now prepared to investigate variational problems and partial differential equations involving space-dependent nonlocal gradients. As an application, we consider a new model in the context  of nonlocal hyperelasticity, characterized by functionals whose stored energy contribution takes the form
	\begin{align}\label{functional_intro}
		\mathcal{F}_{\rho(\cdot)}(u) = \int_{\Omega}f(x, D_{\rho(\cdot)} u)\,dx, 
	\end{align}
where $\Omega\subset \R^n$  is the reference configuration, $u:\Omega\to \R^m$ is an admissible deformation with its heterogeneous nonlocal gradient defined componentwise, and the density $f:\Omega\times \mathbb{R}^{m\times n}\to [0, \infty]$ is a Carath\'{e}odory function with suitable growth and coercivity behavior. 
The deformation maps are assumed to lie in a subspace $X\subset H^{\rho(\cdot), p}(\Omega;\mathbb{R}^m)$, which can encode various choices of 
  boundary conditions. We treat here three different cases, namely, Dirichlet boundary conditions, where 
  $X= H_g^{\rho(\cdot), p}(\Omega;\mathbb{R}^m)$ with given $g\in W^{1-1/p, p}(\partial \Omega;\R^m)$, natural boundary conditions,  i.e., free boundary values, described by $X=\{u\in  H^{\rho(\cdot), p}(\Omega;\mathbb{R}^m): \int_\Omega u\, dx=0\}$, and mixed boundary conditions, where local Dirichlet data is prescribed only on part of the boundary $\partial \Omega$.

Addressing the question of solvability for these vectorial variational problems, we provide sufficient conditions for the existence of minimizers for $\mathcal{F}_{\rho(\cdot)}$  via the direct method.   Indeed, in addition to the coercivity of $\mathcal{F}_{\rho(\cdot)}$, which follows from standard arguments using the nonlocal Poincar\'e inequalities stated in \textit{(h)}, the second key ingredient is the weak lower semicontinuity of these functionals. We show that this property holds when the integrand $f$ is quasiconvex or polyconvex, using \textit{(b)} to reduce the proof to the well-known analogous result in the local setting.  
Since this translation is valid only up to lower-order operators, a careful analysis is needed to handle the mismatch between $\Gr$ and the classical gradient, which can be achieved through strong convergence arguments.
Note also that, as in both the classical and the homogeneous nonlocal settings, quasiconvexity of $f$ remains the natural convexity condition for characterizing the weak lower semicontinuity of $\mathcal{F}_{\rho(\cdot)}$, provided standard $p$-growth from above is assumed for the integrand.

 Finally, we establish a link between the studied nonlocal variational principles and boundary value problems involving heterogeneous nonlocal gradients by deriving optimality conditions in the form of Euler-Lagrange equations.  Our variational approach automatically yields the existence of solutions for these nonlocal systems of partial differential equations, subject to local Dirichlet, Neumann, or mixed boundary conditions.

\subsection{Structure of the manuscript.} The remainder of this paper is organized as follows. Section~\ref{sec:preliminaries} covers the necessary preliminaries, including notational conventions and terminology, as well as a review of the set-up and key results about homogeneous nonlocal gradients. Moreover, it contains relevant elements of Bessel potential spaces and of pseudo-differential operator theory, which, even though well-known, are revisited to ensure a self-contained presentation. 
In Section~\ref{sec:heterogeneousgradients}, we introduce our notion of heterogeneous nonlocal gradient, starting with a definition for smooth functions and extending it to distributions. The latter involves a rigorous characterization in terms of pseudo-differential operators. In this context, we also introduce several useful auxiliary operators, that is, commutators, parametrices, residuals, and extension operators.
As a central part of this work, Section~\ref{sec:heterogeneous_Sobolevspaces} is then devoted to developing a comprehensive theory of the associated heterogeneous nonlocal Sobolev spaces. We provide a detailed discussion and proofs of the tools outlined earlier, see \textit{(a)} - \textit{(h)} above. Finally, Section~\ref{sec:variational} focuses on applications to nonlocal variational problems of the form~\eqref{functional_intro}, where we establish the existence of minimizers for functionals with quasiconvex or polyconvex integrands under local Dirichlet, Neumann, or mixed-type boundary conditions. The corresponding Euler-Lagrange equations are also derived to conclude the analysis.

\section{Preliminaries}\label{sec:preliminaries}

\subsection{General notation and terminology}
In the following, we fix some notations and conventions used throughout the paper.

The Euclidean norm of a vector $x\in \R^n$ is denoted by $|x|$, and $x \cdot y$ is the standard inner product of $x, y\in \R^n$. We write $B_\eps(0)$  for the Euclidean ball of radius $\eps>0$ in $\R^n$ centered at the origin. For a matrix $A\in \R^{N \times n}$,  $|A|$ stands for the Frobenius norm. Let $\langle \xi\rangle:=\sqrt{1+|\xi|^2}$ for $\xi\in \R^n$. It holds that $\max\{1, |\xi|\}\leq \langle\xi \rangle\leq \sqrt{2}\, \max\{1, |\xi|\}$. Moreover, we note for further reference that for $t\geq 0$ and $r \in (0,1]$, 
\begin{equation}\label{eq:abssproperty}
r^t\abss{r\xi}^{-t} \leq \abss{\xi}^{-t} \quad\text{for all $\xi \in \R^n$,}
\end{equation}
which follows from
\begin{align*}
r^t\abss{r\xi}^{-t} = \frac{r^t}{r^t(\frac{1}{r^2}+\abs{\xi}^2)^{t/2}}=\frac{1}{(\frac{1}{r^2}+\abs{\xi}^2)^{t/2}} \leq \abss{\xi}^{-t}.
\end{align*}
For a subset $U \subset \R^n$, its closure is denoted by $\overline{U}$ and $\mathcal{H}^{n-1}(U)$ stands for the $(n-1)$-dimensional Hausdorff measure of $U$. We write $U' \Subset U$ to indicate that $U'\subset \R^n$ is compactly contained in $U$. The identity map on a space $X$ is written as $\Id_X$, or simply $\Id$.  
As common in the literature, we use the symbol $\ast$ for the convolution between two functions. 

Let $X$ be a normed space. 
Regarding the convergence of a sequence $(v_j)_j\subset X$, we write $v_j\to v$ in $X$  and $v_j\weakly v$ in $X$ as $j\to \infty$ to indicate strong and weak convergence, respectively.  The notation $X\hookrightarrow Y$ represents a continuous embedding between two normed spaces $X$ and $Y$, and $X \hookrightarrow\hookrightarrow Y$ a compact embedding.  

Adopting the standard notation for Lebesgue spaces, $L^p(U)$ with $p\in [1, \infty]$ stands for the space of $p$-integrable (essentially bounded for $p=\infty$) functions defined on a measurable set $U\subset \R^n$ with values in $\R$ (or $\C$). The Fourier transform of a (complex-valued) function $f\in L^1(\R^n)$ is defined as
\begin{align*}
\widehat{f}(\xi) = \int_{\R^n} f(x) e^{-2\pi i x\cdot \xi}\, dx  \quad\text{for $\xi\in \R^n$.}
\end{align*}

Moreover, we work with various spaces of smooth functions. These include the Schwartz space $\Scal(\R^n)$ and for an open set $U\subset \R^n$, the space $C^\infty_c(U)$ of smooth test functions with compact support in $U$.  Moreover, by $C^{\infty}(\overline{U})$ we denote the restriction of all smooth functions on $\R^n$ to $U$, and $C^{\infty}_0(\overline{U})$ denotes all functions in $C^{\infty}(\overline{U})$ whose extension to $\R^n$ by zero is smooth.

For the partial derivatives of smooth function $f: U\to \R, \ x\mapsto f(x)$, we use the standard multi-index notation $\partial^\alpha f = \partial_{x_1}^{\alpha_1}\partial_{x_2}^{\alpha_2}\ldots \partial_{x_n}^{\alpha_n} f$ for $\alpha\in \N_0^n$. The gradient of $f$ is denoted by $\nabla f=(\partial_{x_1} f, \ldots, \partial_{x_n} f)$. For a function $f$ of two vector variables $x\in \R^n$ and $\xi\in \R^n$, we write $\nabla_x f$ and $\nabla_\xi f$ for the collections of partial derivatives $(\partial_{x_1}f, \ldots, \partial_{x_n}f)$ and $(\partial_{\xi_1}f, \ldots, \partial_{\xi_n} f)$, respectively. The same notations are used for distributional and weak derivatives when $f$ is less regular.

We denote by $\Dcal'(U)$ the space of distributions on an open set $U\subset \R^n$, defined as the dual space of $C_c^\infty(U)$, and by   $\mathcal{S}'(\R^n)$ the space of tempered distributions, i.e., the dual space of $\Scal(\R^n)$.

For Sobolev spaces, we also use the common notation, that is, for $p\in [1, \infty]$ and $U\subset \R^n$ open, the space $W^{1,p}(U)$ comprises all functions $f\in L^p(U)$ whose distributional gradient $\nabla f$ lies in $L^p(U;\R^n)$, endowed with the norm  $\norm{u}_{W^{1,p}(U)}=\norm{u}_{L^p(U)} + \norm{\nabla u}_{L^p(U;\R^{n})}$. If $\Omega\subset \R^n$ is a bounded Lipschitz domain and $p\in (1, \infty)$, then $W^{1,p}(\Omega)$ admits an associated trace space, namely the Sobolev-Slobodeckij space $W^{1-1/p,p}(\partial \Omega)$. We denote by $W^{1,p}_g(\Omega)$ the subspace of $W^{1,p}(\Omega)$ consisting of all functions with prescribed trace $g\in W^{1-1/p,p}(\partial \Omega)$. 

Without explicit mention, functions $f: U\subset \R^n\to \R$ are identified with their extension by zero outside $U$, where needed. For instance, $C_c^\infty(U)$ and $C^{\infty}_0(\overline{U})$ with $U\subset \R^n$ open are typically considered a subspace of $C^\infty(\R^n)$ or $C_c^\infty(\R^n)$; to emphasize this identification, we occasionally write $\tilde f=\mathbbm{1}_U f$ for the trivial extension of $f$.

The pointwise restriction of a function $\tilde f:\R^n\to \R$ to a set $U\subset \R^n$ is denoted by $\tilde f|_U$.  If $X$ is a function space on $\R^n$,  then $X|_{U}$ is the space of restrictions $\tilde f|_U$ for $\tilde f\in X$. Where necessary, we implicitly identify functions in $X$ with their restrictions to $U\subset \R^n$ without further comment. 
For the restriction of a distribution $\tilde f\in \D'(\R^n)$ to an open set $U\subset \R^n$, we define $\tilde f|_U\in \Dcal'(U)$ via 
\begin{align*}
 \langle \tilde f|_U ,  \varphi \rangle = \langle \tilde f,  \varphi\rangle \quad \text{for all $\phi\in C^\infty_c(U)\subset C_c^\infty(\R^n)$.}
\end{align*} 
Note that the restriction operator 
\begin{align}\label{eq:piOmega}
\pi_\Omega: \Dcal'(\R^n)\to \Dcal'(U),  \ \pi_\Omega f:=\tilde f|_U \quad \text{for $\tilde f\in \Dcal'(\R^n)$, }
\end{align} is linear and continuous. When applied to Lebesgue functions, $\pi_\Omega: L^p(\R^n)\to L^p(U)$ is also bounded and corresponds to  the pointwise restriction to $U$.

All the above-mentioned properties of functions and the definitions of function spaces naturally extend to vector-valued functions by applying them componentwise.

For real functions, we understand the notions of increasing and decreasing in the non-strict sense.  
A function $f:\R\to \R$ is called almost decreasing if there is a $C>0$ such that $f(t) \geq Cf(s)$ for $t \leq s$, and an analogous definition holds for almost increasing.

Finally, unless mentioned otherwise, $C, c$ are positive constants that may vary from line to line. When the dependence on a specific parameter is relevant, it will be indicated with a subscript; for example, $c_\alpha$ or $ C_{\alpha}$ for a constant depending on $\alpha$. 

\subsection{Nonlocal gradients}\label{sec:nonlocalgradients}
Let us recall the setting of general nonlocal gradients with radial kernels as it was introduced in~\cite{BeMCSc24}, and state some of their properties. This will later serve as the basis for defining nonlocal gradients with heterogeneous horizons.

Throughout the paper, take $\rho_1:\R^n \setminus\{0\} \to [0,\infty)$ to be a radial and integrable function with compact support in $B_1(0)$, satisfying
\begin{enumerate}[label = (H\arabic*)]
\setcounter{enumi}{-1}
\item \label{itm:h0} $\inf_{\overline{B_\epsilon(0)}} \rho_1 >0$ for some $\eps >0$ and $\displaystyle \int_{\R^n} \rho_1 \,dz =n$.
\end{enumerate}

 We denote by $\overline{\rho}_1:(0,\infty) \to [0,\infty)$ its radial representative, that is, $\rho_1(x)=\overline{\rho}_1(\abs{x})$ for all $x \in \R^n \setminus \{0\}$. Following~ \cite{BeMCSc24} (cf. also~\cite{CKS25}), we further assume the conditions \ref{itm:h1}-\ref{itm:upper}: \smallskip
\begin{enumerate}[label = (H\arabic*)]

\item \label{itm:h1} The function $f_{\rho_1}:(0,\infty)\to \R, \ r \mapsto r^{n-2}\overline{\rho}_1(r)$ is decreasing on $(0,\infty)$ and $r \mapsto r^\iota  f_{\rho_1}(r)$ is decreasing on $(0,\epsilon)$ for some $\iota>0$;

\item \label{itm:derivatives} $f_{\rho_1}$ is smooth on $(0,\infty)$ and for every $k\in \N$ there exists a $C_k>0$ with 
\[
\abslr{\frac{d^k}{dr^k}f_{\rho_1}(r)} \leq C_k \frac{f_{\rho_1}(r)}{r^k} \quad \text{for $r\in (0, \epsilon)$};
\]

\item \label{itm:lower} the function $r \mapsto r^{n+\lambda-1}\overline{\rho}_1(r)$ is almost decreasing on $(0,\epsilon)$; 

\item \label{itm:upper} the function $r \mapsto r^{n+\kappa-1}\overline{\rho}_1(r)$ is almost increasing on $(0,\epsilon)$.
\end{enumerate} \smallskip

Roughly speaking, the hypotheses \ref{itm:lower} and~\ref{itm:upper} mean that $\rho_1$ lies in between fractional kernels of orders $\lambda$ and $\kappa$ with $\lambda\leq \kappa$, close to the origin.

Here are two examples of kernels that fulfill these assumptions (see~\cite{BeMCSc24} for more) and will be used to illustrate our results later.
\begin{example}\label{ex:rho1} Let $s\in (0,1)$ and $w \in C_c^{\infty}(B_1(0))$ be a non-zero radial function such that $w(0)>0$ and $w/\abs{\cdot}^{1+s}$ is radially decreasing.
Then, the kernels $\rho_1$ given by
\[
\rho_1(z) = \frac{w(z)}{\abs{z}^{n+s-1}}\quad  \text{and}\quad \rho_1(z)=\frac{w(z)\log(1/\abs{z})}{\abs{z}^{n+s-1}}\qquad \text{for $z\in \R^n\setminus\{0\}$}
\]
meet~\ref{itm:h1}-\ref{itm:upper} with parameters $\iota=n-2$, $\lambda=s$, and  $\kappa=s$  in the first case, and any $\kappa \in (s,1)$ in the second case.
After rescaling by a suitable constant, both kernels satisfy also the hypothesis~\ref{itm:h0}.
\end{example}

The nonlocal gradient associated to the kernel function $\rho_1$ is defined as follows.
\begin{definition}[Nonlocal gradient]\label{def:nonlocalgradient}
For $\phi \in \Scal(\R^n)$, the nonlocal gradient $D_{\rho_1} \varphi:\R^n \to \R^n$ is given by
\[
D_{\rho_1}\phi(x) = \int_{\R^n} \frac{\phi(x)-\phi(y)}{\abs{x-y}}\frac{x-y}{\abs{x-y}}\rho_1(x-y)\,dy,  \quad x \in \R^n.
\]
\end{definition}
\color{black}
There is an alternative definition of the nonlocal gradient via convolution. Indeed, with the function
\begin{equation}\label{eq:Qrhodefinition}
Q_{\rho_1}:\R^n \setminus \{0\} \to [0,\infty), \quad z\mapsto \int_{\abs{z}}^1 \frac{\overline{\rho}_1(r)}{r}\,dr,
\end{equation}
which is integrable with $\norm{Q_{\rho_1}}_{L^1(\R^n)}=1$ and has compact support in $B_1(0)$, it holds for all $\phi \in \Scal(\R^n)$ that
\begin{equation}\label{eq:translation}
D_{\rho_1}\phi = Q_{\rho_1} * \nabla \phi = \nabla (Q_{\rho_1} * \phi),
\end{equation}
cf.~\cite[Lemma~2.5, Proposition~2.6]{BeMCSc24}. As a consequence, one obtains the Fourier representation
\begin{equation*}\label{eq:fourierrepr}
\widehat{D_{\rho_1} \phi}(\xi) =  \widehat{Q}_{\rho_1}(\xi) 2\pi i\xi \widehat{\phi}(\xi) \quad \text{for all $\xi \in \R^n$},
\end{equation*}\color{black}
noting that $D_{\rho_1} \phi \in \Scal(\R^n;\R^n)$.
The following bounds on the Fourier transform $\widehat{Q}_{\rho_1}$ and its derivatives, which result from~\cite[Lemma~4.3, Lemma~4.10, Lemma~7.1]{BeMCSc24},  will play a crucial role in developing our theory. 

\begin{proposition}\label{prop:BeMCSc24}
Let $\rho_1$ satisfy \ref{itm:h0}-\ref{itm:upper}. Then $\widehat{Q}_{\rho_1}$ is smooth, positive, and there is a $C>0$ such that
\[
\frac{1}{C} \frac{\bar{\rho}_1(1/\abs{\xi})}{\abs{\xi}^n} \leq \widehat{Q}_{\rho_1}(\xi) \leq C \frac{\bar{\rho}_1(1/\abs{\xi})}{\abs{\xi}^n} \quad \text{for all $\xi\in \R^n$ with $\abs{\xi} \geq 1/\epsilon$}.
\]
Moreover, for every $\alpha \in \N_0^n$, it holds that
\begin{align}\label{est2:lemmaBeMCSc24}
\absb{\partial^{\alpha} \widehat{Q}_{\rho_1}(\xi)} \leq C_{\alpha} \abss{\xi}^{-\abs{\alpha}} \absb{\widehat{Q}_{\rho_1}(\xi)} \quad \text{for all $\xi \in \R^n$}.
\end{align}
\end{proposition}

\subsection{Bessel potential spaces}\label{sec:bps}
In this section, we collect the necessary background on Bessel potential spaces on $\R^n$, as well as on bounded domains.
 These spaces naturally arise in the theory of pseudo-differential operators.

By definition, the Bessel potential space for $s \in \R$ and $p \in (1,\infty)$ is
\[
H^{s,p}(\R^n) =\{ u \in \Scal'(\R^n) \,:\, (\abss{\cdot}^s \widehat{u})^{\vee} \in L^p(\R^n)\},
\]
with the norm $\norm{u}_{H^{s,p}(\R^n)}=\norm{(\abss{\cdot}^s \widehat{u})^{\vee}}_{L^p(\R^n)}$. For more on these spaces, we refer~e.g.~ to~\cite{Trie1, Gra14b}, noting that they fit into the more general setting of Triebel-Lizorkin spaces, that is, $H^{s,p}(\R^n)=F^{s}_{p,2}(\R^n)$. If $s<0$, then $H^{s,p}(\R^n)$ is a space of tempered distributions, while for $s \geq 0$, its elements can be identified with functions in $L^p(\R^n)$; when $s$ is a positive integer, the spaces coincide with the classical Sobolev spaces, that is, $H^{k, p}(\R^n)= W^{k,p}(\R^n)$ for $k\in \N$. It is worth mentioning here that an equivalent norm on $H^{s, p}(\R^n)$ is given by
\begin{align}\label{eq:eqnormrn}
\norm{u}_{H^{s,p}(\R^n)} \cong \norm{u}_{H^{s-1,p}(\R^n)}+\norm{\nabla u}_{H^{s-1,p}(\R^n;\R^n)},
\end{align}
see~\cite{Trie1}.

For an open and bounded set $U \subset \R^n$,  one can also define the restricted space
\begin{align*}
H^{s,p}(U)=\{u \in \Dcal'(U) \,:\, u = \tilde{u}|_{U} \ \text{for some $\tilde{u} \in H^{s,p}(\R^n)$}\} 
\end{align*}
equipped with the norm
\begin{align}\label{normHsp}
\norm{u}_{H^{s,p}(U)}=\inf \left\{\norm{\tilde{u}}_{H^{s,p}(\R^n)} \,:\, \tilde{u} \in H^{s,p}(\R^n), \ \tilde{u}|_{U}=u\right\}.
\end{align}
If $s\geq 0$, then $H^{s,p}(U) = H^{s, p}(\R^n)|_{U}\subset L^p(U)$ and we also define  the Bessel potential spaces with zero complementary values as 
\begin{align}\label{Besselzero}
H^{s,p}_0(U) = \{ u \in H^{s,p}(\Rn)\,:\, u = 0 \ \text{a.e.~in $U^c$}\};
\end{align}
with our notational convention about trivial extensions, $H^{s, p}_0(U)$ can be seen as both a subspace of $H^{s,p}(\R^n)$ and $H^{s,p}(U)$. Note that, by density, $H_0^{s,p}(U)$ can also be identified with the closure of $C_c^\infty(U)$ in $H^{s,p}(\R^n)$ if $U$ is a Lipschitz domain, see~e.g.~\cite[Theorem~3.9\,(iii) and Example~3.11\,(b)]{BeMCSc24}.

A key property that we will often use is that $H^{s_1,p}(U)$ embeds compactly into $H^{s_2,p}(U)$ for any $s_1, s_2\in \R$ with $s_1>s_2$,  in formulas,
\begin{align}\label{Hcompactembedding}
H^{s_1, p}(U) \hookrightarrow\hookrightarrow H^{s_2, p}(U) \qquad \text{for $s_1>s_2$}; 
\end{align}
this should be a well-known fact, but for general open and bounded sets  we could only find it as a special case of  the more general setting of~\cite[Theorem~3.2]{GHS21}.
Moreover, analogously to~\eqref{eq:eqnormrn}, the equivalence of norms 
\begin{align}\label{eq:besselboundednorm}
\norm{u}_{H^{s,p}(\Omega)} \cong \norm{u}_{H^{s-1,p}(\Omega)}+\norm{\nabla u}_{H^{s-1,p}(\Omega;\R^n)}
\end{align}
holds if $\Omega$ is a bounded Lipschitz domain, see~\cite[Theorem~1.1]{ShY23}.

An important tool for our analysis is an extension operator connecting the Bessel potential spaces on $\Omega$ and $\mathbb{R}^n$. Precisely, for a bounded Lipschitz domain $\Omega \subset \mathbb{R}^n$, there exists a linear extension $E$ such that
\begin{align}\label{Rychkov}
\text{$E: H^{s, p}(\Omega)\to H^{s,p}(\R^n)$ is bounded for all $s\in \R$ and all $p\in(1,\infty)$.}
\end{align} 
This result is a consequence of~\cite[Theorem~4.1]{Ryc99} by Rychkov, where a universal linear extension operator  bounded on all Triebel-Lizorkin spaces is established.

Let us conclude this subsection by introducing some shorthand notations for bounded linear operators mapping between Bessel potential spaces. We say that an operator $\mathcal{L}$ lies in the class $\LBH_\Omega^\gamma$ with $\gamma\in \R$, if 
\begin{align*}
\text{$\Lcal:H^{s,p}(\Omega) \to H^{s-\gamma,p}(\Omega)$ \quad  is bounded for all $s \in \R$ and $p \in (1,\infty)$; }
\end{align*}
further, we set
$\LBH_\Omega^{-\infty} := \bigcap_{\gamma\in \R} \LBH_\Omega^\gamma$ and $\LBH_\Omega^{+} := \bigcap_{\gamma>0} \LBH_\Omega^\gamma$.  
We use the analogous notations with $\R^n$ in place of $\Omega$, dropping the subscript, e.g., $\LBH^\gamma$ consists of all operators $H^{s,p}(\R^n) \to H^{s-\gamma,p}(\R^n)$ that are bounded for all $s \in \R$ and $p \in (1,\infty)$.

\subsection{Pseudo-differential operators}\label{sec:pdo}
We next recall the definition of pseudo-differential operators with symbols in the H\"ormander class, along with some of their key properties. These operators are fundamental to the study of heterogeneous nonlocal gradients, as the latter can be shown to fit into this abstract framework (see Lemma~\ref{le:symbolbounds}). Our focus here lies on gathering the results necessary to ensure a self-contained presentation. For a broader and more detailed introduction to the topic, we refer the reader to standard textbooks such as~\cite{Kum81, Hor07}.

Throughout this section, all functions are considered complex-valued. The H\"{o}rmander class of symbols $S^{m}_{\nu,\mu}$ with $m \in \R$, $0\leq \mu < \nu \leq 1$ consists of smooth functions $a \in C^{\infty}(\R^n\times\R^n)$ such that for all $\alpha,\beta \in \N_0^n$
\[
\abs{\partial_x^\beta \partial_\xi^\alpha a(x,\xi)} \leq C_{\alpha,\beta}\abss{\xi}^{m-\nu\abs{\alpha}+\mu\abs{\beta}} \quad\text{for all $x,\xi \in \R^n$.}
\]
A family of seminorms on $S^{m}_{\nu,\mu}$ given by
\[
\abs{a}_{S^{m}_{\nu,\mu},k}:=\sup_{\abs{\alpha+\beta}\leq k} \sup_{x,\xi \in \R^n} \abs{\partial_x^\beta \partial_\xi^\alpha a(x,\xi)}\abss{\xi}^{-m+\nu\abs{\alpha}-\mu\abs{\beta}}
\]
for $k\in \N_0$, which turns $\Smm$ into a Fr\'{e}chet space.

Any symbol $a \in \Smm$ induces a pseudo-differential operator  $\Acal=\Op(a)$ that is continuous from $\Scal(\R^n)$ to $\Scal(\R^n)$ and defined by
\begin{align*}
\Acal \phi (x) = \int_{\R^n} e^{2\pi i x \cdot \xi}a(x,\xi)\widehat{\phi}(\xi)\,d\xi = \int_{\R^n}\int_{\R^n}e^{2\pi i(x-y)\cdot\xi}a(x,\xi)\phi(y)\,dy\,d\xi \quad \text{for $x \in \R^n$}.
\end{align*}
Conversely, for any such operator $\Acal$, we write $\Acal \in \OPS^{m}_{\nu,\mu}$ and denote its symbol, which is uniquely defined by the operator, by $\sigma(\Acal)$.

For $\Acal \in \OPSmm$, there always exists an adjoint operator $\mathcal{A}^* \in \OPSmm$ such that
\[
\int_{\R^n} \Acal \phi \, \psi^*\,dx = \int_{\R^n} \phi \, (\Acal^* \psi)^*\,dx \quad \text{for all $\phi,\psi \in \Scal(\R^n)$},
\]
with $(\cdot)^*$ denoting complex conjugation, see~e.g.~\cite[Theorem~2.1.7]{Kum81} also for an asymptotic formula of $\sigma(\Acal^*)$. This enables $\Acal$ to be extended to a continuous map from $\Scal'(\R^n)$ to $\Scal'(\R^n)$ via duality as
\[
\inner{\Acal u , \phi} := \inner{u, \Acal^* \phi} \quad \text{for all $\phi \in \Scal(\R^n)$},
\]
with $u \in \Scal'(\R^n)$.  More sophisticated boundedness properties of pseudo-differential operators on the Bessel potential spaces hold as well. For our purposes, we observe that any $\Acal \in \OPS^{m}_{1,\mu}$ satisfies that 
\begin{align*}
\text{$\Acal:H^{s,p}(\R^n) \to H^{s-m,p}(\R^n)$ is bounded for all $s \in \R$ and $p \in (1,\infty)$, }
\end{align*}
see~e.g.~\cite[Proposition~6.5]{Tay11b}; using our notation from the previous section, this means $\Acal\in \LBH^m$. In fact, for every $s\in \R$ and $p \in (1,\infty)$ there is some $k \in \N$ and a constant $C>0$ such that
\begin{equation}\label{eq:lpboundseminorm}
\norm{\Acal u}_{H^{s-m,p}(\R^n)} \leq C\abs{\sigma(\Acal)}_{S^{m}_{1,\mu},k}\norm{u}_{H^{s,p}(\R^n)} \quad \text{for all $u \in H^{s,p}(\R^n)$},
\end{equation}
cf.~\cite{Fef73,Miy87} or the abstract argument in \cite[Lemma~3.5]{AbP18}. \smallskip

The following paragraphs discuss three selected topics from the theory of pseudo-differential operators, tailored to meet the requirements of our analysis.\smallskip

\textit{a) Composition of symbols.} For two symbols $a_j \in S^{m_j}_{\nu,\mu}$ for $j=1,2$, there exists a symbol $a_1\circ a_2 \in S^{m_1+m_2}_{\nu,\mu}$ such that $\Op(a_1)\Op(a_2) = \Op(a_1\circ a_2)$. In fact, by \cite[Chapter 2, Theorem~3.1]{Kum81} with $N=1$, the following formula holds for the symbol of the composition: For all $x, \xi\in \R^n$,
\begin{align*}
(a_1\circ a_2)(x,\xi) = a_1(x,\xi)a_2(x,\xi) + \sum_{l=1}^n \int_0^1 r_{l,\theta}(x,\xi)\,d\theta,
\end{align*}
 with $r_{l,\theta}$ given by
\[
r_{l,\theta}(x,\xi)=\mathrm{Os}\text{-}\int_{\R^{n}\times \R^n} e^{-2\pi i y \cdot \eta} \partial_{\xi_l}a_1(x,\xi+\theta\eta)\partial_{x_l} a_2(x+y,\xi)\,d(y,\eta),
\]
which is an oscillatory integral (see e.g.~\cite[Chapter 1, \S 6]{Kum81}); however, we will not go into the details of its definition. More importantly,~\cite[Chapter 2, Lemma~2.4]{Kum81} implies that for any $k \in \N$ there exist a constant $C>0$ and an integer $k'\in \N$, both independent of $\theta$ and $l$, such that
\[
\abs{r_{l,\theta}}_{S^{m_1+m_2-\nu+\mu}_{\nu,\mu},k}\leq C\abs{\partial_{\xi_l}a_1}_{S^{m_1-\nu}_{\nu,\mu},k'}\abs{\partial_{x_l}a_2}_{S^{m_2+\mu}_{\nu,\mu},k'};
\] 
hence, by integrating over $\theta$, it follows that $a_1\circ a_2-a_1a_2 \in S^{m_1+m_2-\nu+\mu}_{\nu,\mu}$ with
\begin{align}\label{eq:compositionbound}
\abs{a_1\circ a_2-a_1a_2}_{S^{m_1+m_2-\nu+\mu}_{\nu,\mu},k}\leq C\abs{\nabla_\xi a_1}_{S^{m_1-\nu}_{\nu,\mu},k'}\abs{\nabla_x a_2}_{S^{m_2+\mu}_{\nu,\mu},k'}.
\end{align}
\smallskip

\color{black}
\textit{b) Parametrices and invertibility.} Given a symbol $a \in S^{m}_{\nu,\mu}$, it is well-known that one can construct under suitable assumptions a so-called parametrix for $\Op(a)$, which is an almost inverse. 
 
Suppose that $a$ satisfies the following conditions:
There exists $m'\in \R$ with $m' \leq m$ and a constant $C>0$ such that
\begin{equation}\label{eq:parametrixcond1}
\abs{a(x,\xi)} \geq C\abss{\xi}^{m'} \quad \text{for all $x ,\xi \in \R^n$,}
\end{equation}
and for all multi-indices $\alpha,\beta \in \N_0^n$ there are constants $C_{\alpha, \beta}>0$ such that
\begin{equation}\label{eq:parametrixcond2}
\abs{\partial_x^\beta\partial_\xi^\alpha a(x,\xi)} \leq C_{\alpha,\beta}\abss{\xi}^{-\nu\abs{\alpha}+\mu\abs{\beta}}\abs{a(x,\xi)}  \quad \text{for all $x ,\xi \in \R^n$.}
\end{equation}
 Then, there exists a parametric symbol $b \in S^{-m'}_{\nu,\mu}$ such that
\begin{align*}
\Op(b) \Op(a)-\Id=\Op(r) \quad \text{and}\quad \Op(b)\Op(a)-\Id=\Op(r'),
\end{align*}
where $r,r' \in S^{-\infty}:=\bigcap_{m \in \R} S^{m}_{1,0}$ are so-called smoothing symbols, see e.g.~\cite[Chapter~2, Theorem~5.4]{Kum81}, and $\text{Id}$ is the identity operator, which can be viewed as pseudo-differential operator associated to the constant function with value $1$. We write $\Op(r),\Op(r') \in \OPS^{-\infty}$. Note that these operators are bounded from $H^{s,p}(\R^n)$ to $H^{s-\gamma,p}(\R^n)$ for any $s,\gamma\in \R$ and $p \in (1,\infty)$, that is,
\begin{align*}
\Op(r),\Op(r')\in \LBH^{-\infty}.
\end{align*} 

The exact invertibility of pseudo-differential operators is much harder in general, but there is a suitable replacement for slowly varying symbols.  This is the following statement, which is an adaptation of \cite[Theorem~14 and 15]{RaR07} to the symbol classes that are relevant for us. 
\begin{proposition}\label{prop:exactinvertible}
Let $a \in S^{0}_{1,\mu}$ satisfy \eqref{eq:parametrixcond1} with $m'\geq \mu-1$ and~ \eqref{eq:parametrixcond2}, and let $a_\epsilon:=a(\cdot,\epsilon\,\cdot) \in S^{0}_{1,\mu}$ for $\epsilon >0$. Then, for any $p\in (1, \infty)$ there exists a $\eps_0>0$ such that for every $\eps\in (0, \eps_0)$ it holds that
\[
\Op(a_\epsilon^{-1})\Op(a_\epsilon):L^p(\R^n) \to L^p(\R^n) \ \text{is invertible}.
\]
In particular, $\Op(a_\eps):L^p(\R^n)\to L^p(\R^n)$ has a trivial kernel. 
\end{proposition}

\begin{proof}
One may assume without loss of generality that $m' <0$, since \eqref{eq:parametrixcond1} will still be satisfied. Note that by virtue of \eqref{eq:parametrixcond1} and~ \eqref{eq:parametrixcond2}, it holds that $a^{-1} \in S^{-m'}_{1,\mu}$, and thus also $a_\eps^{-1}\in S^{-m'}_{1,\mu}$ for every $\eps>0$.

Let $k\in \N_0$. We compute for all multi-indices $\alpha,\beta \in \N_0^n$ with $\abs{\alpha+\beta}\leq k$ that
\begin{align*}
\abs{\partial^\beta_x \partial^\alpha_\xi a_\epsilon(x,\xi)} &= \epsilon^{\abs{\alpha}}\abs{\partial^\beta_x \partial^\alpha_\xi a(x,\epsilon \xi)} \\
&\leq \abs{a}_{S^{0}_{1,\mu},k}\epsilon^{\abs{\alpha}}\abss{\epsilon \xi}^{-\abs{\alpha}+\mu\abs{\beta}} \\ &\leq \abs{a}_{S^{0}_{1,\mu},k} \abss{\xi}^{-\abs{\alpha}}\abss{\xi}^{\mu\abs{\beta}} = \abs{a}_{S^{0}_{1,\mu},k}\abss{\xi}^{-\abs{\alpha}+\mu\abs{\beta}} 
\end{align*}
for $x, \xi\in \R^n$; here, the last line uses \eqref{eq:abssproperty}  as well as $\mu\abs{\beta}\geq 0$, assuming that $\epsilon <1$. A similar estimate with $\alpha \not =0$ shows that
\begin{align*}
\abs{\partial^\beta_x \partial^\alpha_\xi (a_\epsilon)^{-1}(x,\xi)} &\leq \abs{\nabla_\xi a^{-1}}_{S^{-m'-1}_{1,\mu},k}\epsilon^{\abs{\alpha}} \abss{\epsilon\xi}^{-m'-\abs{\alpha}+\mu\abs{\beta}}\\
&\leq \abs{\nabla_\xi a^{-1}}_{S^{-m'-1}_{1,\mu},k}\epsilon^{-m'} \abss{\xi}^{-m'-\abs{\alpha}+\mu\abs{\beta}},
\end{align*}
where~\eqref{eq:abssproperty} has been applied for $t=\abs{\alpha}+m'\geq \mu>0$. 
As a consequence of these computation, it follows that 
\begin{align*}
\abs{\nabla_x a_\epsilon}_{S^{\mu}_{1,\mu},k} \leq \abs{\nabla_x a}_{S^{\mu}_{1,\mu},k} 
\quad \text{and}  \quad \abs{\nabla_\xi (a_\epsilon)^{-1}}_{S^{-m'-1}_{1,\mu},k} \leq \epsilon^{-m'}\abs{\nabla_\xi a^{-1}}_{S^{-m'-1}_{1,\mu},k}. 
\end{align*}

We now obtain along with~\eqref{eq:compositionbound} that
\begin{align*}
 \abs{(a_\epsilon)^{-1} \circ a_\epsilon - 1}_{S^{0}_{1,\mu},k} &\leq  \abs{(a_\epsilon)^{-1} \circ a_\epsilon - 1}_{S^{-m'-1+\mu}_{1,\mu},k}  \leq C \eps^{-m'} \abs{\nabla_\xi a^{-1}}_{S^{-m'-1}_{1,\mu},k'}\abs{\nabla_x a}_{S^{\mu}_{1,\mu},k'}, 
\end{align*}
for some $k'\in \N$ and a constant $C>0$, where the first inequality uses that $-m'-1+\mu \leq 0$.

In view of \eqref{eq:lpboundseminorm}, this shows that $\Op({a_\epsilon^{-1}})\Op({a_\epsilon})-\Id$ is a bounded operator from $L^{p}(\R^n)$ to $L^{p}(\R^n)$ for any $p\in (1, \infty)$, whose operator norm bounded by a constant times $\epsilon^{-m'}$.  
By choosing $\epsilon$ small enough such that this operator norm is smaller than $1$, we conclude the statement via a Neumann series argument. 
\end{proof}

\color{black}

\textit{c) Restriction of pseudo-differential operators.}
This paragraph introduces operators on function spaces over bounded domains in $\R^n$ by suitably restricting pseudo-differential operators, an approach that will be crucial for our analysis later on. 

Let $\Omega\subset \R^n$ be an open, bounded set
and consider a pseudo-differential operator $\Acal \in \OPSmm$. 
If the dual $\Acal^*\in \OPSmm$ maps $C_c^{\infty}(\Omega)$ 
to itself continuously, then we can define $\Acal$ as an operator from $\Dcal'(\Omega)$ 
to itself via duality, denoted by $\Acal^{\O}$; that is, for $u\in \Dcal'(\Omega)$, 
\begin{align}\label{eqrestrict}
\langle \Acal^\Omega u, \psi\rangle = \langle u, \Acal^* \psi \rangle \quad \text{for all $\psi\in C_c^\infty(\Omega)$. }
\end{align}

If $\Acal \in \textrm{OPS}^{m}_{1,\mu}$, then one has that $\Acal^\Omega\in \LBH_\Omega^m$, meaning
\begin{align}\label{mappingproperties_AcalOmega}
\text{$\Acal^{\Omega}:H^{s,p}(\Omega) \to H^{s-m,p}(\Omega)$ is bounded for all $s \in \R$ and $p \in (1,\infty)$. }
\end{align}
To see this, let $u \in H^{s,p}(\Omega)$. Observe that if $\tilde{u} \in H^{s,p}(\R^n)$ is such that $\tilde{u}|_{\Omega}=\pi_\Omega \tilde u=u$, then $(\Acal \tilde{u})|_{\Omega}$ agrees with $\Acal^\O u$. In view of~\eqref{normHsp}, we conclude along with~\eqref{eq:lpboundseminorm} that
\begin{align*}
\norm{\Acal^\O u}_{H^{s-m,p}(\Omega)} &= \inf\{\norm{\tilde v}_{H^{s-m, p}(\Omega)}\,:\,\tilde v\in H^{s-m, p}(\R^n), \tilde v|_\Omega=\Acal^\O u\} \\ 
& \leq \inf\{\norm{\Acal \tilde{u}}_{H^{s-m,p}(\R^n)}\,:\, \tilde{u} \in H^{s,p}(\R^n), \tilde{u}|_{\Omega}=u\}
 \\ & \leq C \inf\{\norm{\tilde{u}}_{H^{s,p}(\R^n)} \,:\, \tilde{u} \in H^{s,p}(\R^n), \tilde{u}|_{\Omega}=u\}= \norm{u}_{H^{s,p}(\Omega)}
\end{align*}
with a constant $C>0$, which shows the claim.

\subsection{Construction of smooth horizon function}
We will establish here the existence of suitable smooth horizon functions that satisfy the requirements important for defining our heterogeneous nonlocal gradients with desirable properties later on. In particular, these functions must decay sufficiently fast, essentially at an exponential rate, toward the boundary.

\color{black}
\begin{lemma}\label{le:constructioneta}
Let $\Omega \subset \R^n$ be a bounded Lipschitz domain. There exists a non-negative function $ \delta \in \Cinfo$ such that the following properties hold:
\begin{enumerate}[label=(\roman*)]
\item $\delta >0$ on $\Omega$,\\[-0.3cm]
\item $\delta(x) \leq \dist(x,\Omega^c)$ for all $x \in \R^n$, \\[-0.3cm]
\item for all $\gamma\in (0,1)$ and every $\alpha \in \N_0^n$  there is a constant $C>0$ such that
\begin{align}\label{eq:decayeta}
\abs{\partial^{\alpha} \delta(x)} \leq C \abs{\delta(x)}^{1-\gamma} \quad \text{for all $x \in \R^n$}.
\end{align}
\end{enumerate}
\end{lemma}

\begin{proof}
Let $d \in C^{\infty}(\Omega)$ satisfy $c_1 \dist(x,\Omega^c) \leq d(x) \leq c_2 \dist(x,\Omega^c)$ with $0<c_1\leq c_2$ for all $x \in \Omega$, and also
\begin{align}\label{eq:distancefunction}
\abs{\partial^{\alpha} d(x)}\leq C \abs{d(x)}^{1-\abs{\alpha}} \quad \text{ for all $\alpha \in \N_0^n$ and $x \in \Omega$,}
\end{align}
see e.g.,~\cite[Theorem VI.2]{Stein} for a construction. Let us define 
\[
\delta(x)=\begin{cases}
c e^{-1/d(x)} &\text{if $x \in \Omega$,}\\
0 &\text{if $x \in \Omega^c$},
\end{cases}
\]
with $c>0$ chosen such that $\delta(x) \leq \dist(x,\Omega^c)$ for all $x \in \R^n$. By the Leibniz rule and \eqref{eq:distancefunction}, one can then find  for each $\alpha \in \N_0^n$ a $k \in \N_0$ such that
\begin{align*}
\abs{\partial^{\alpha} \delta(x)} \leq C \abs{d(x)}^{-k} e^{-1/d(x)} = C \abs{d(x)}^{-k}e^{-\gamma/d(x)} \abs{\delta(x)}^{1-\gamma} \quad \text{for all $x \in \Omega$}.
\end{align*}
Since the function $\abs{d}^{-k}e^{-\gamma/d}$ is bounded on $\Omega$ for all $\gamma>0$, we find that \eqref{eq:decayeta} follows. The other properties are clear.
\end{proof}

The proof gives indication on how to construct specific horizon functions, a basic example is the following.
\begin{example} Suppose $\Omega=B_1(0)$. 
A simple admissible choice of horizon function is for instance 
\begin{center}
$\delta(x)=e^{\frac{1}{|x|^2-1}}$\quad  for $x\in B_1(0)$. 
\end{center}
Also interesting from a modeling perspective are horizon functions that are constant in the interior and decay to zero in a boundary layer, which can be obtained via smooth interpolation. 
\end{example}

\section{Heterogeneous horizon nonlocal gradients}\label{sec:heterogeneousgradients}

We start this section by introducing the heterogeneous nonlocal gradients, and show that they fit into the theory of pseudo-differential operators with symbols in the H\"{o}rmander class. This enables us to understand some mapping properties of these operators and to prove an integration by parts formula with local boundary conditions.

\subsection{Definition and first properties of the heterogeneous nonlocal gradient}\label{subsec:defnonlocalgradient}
Throughout the paper, $\Omega \subset \R^n$ is a bounded  and connected Lipschitz domain,  $p\in (1, \infty)$ and $\rho_1$ is a kernel as in Section~\ref{sec:nonlocalgradients}, satisfying the hypotheses~\ref{itm:h0} and~\ref{itm:h1}-\ref{itm:upper}. Further, we introduce a position-dependent horizon function $ \delta \in \Cinfo$ as in Lemma~\ref{le:constructioneta}, which fulfills the assumptions $(i)$-$(iii)$; that is, $\d$ is positive on $\Omega$, bounded from above by the distance to $\Omega^c$, and 
for every $\alpha \in \N_0^n$ and $\gamma\in (0,1)$,  there is a constant $C>0$ such that
\begin{align}\label{eq:decaybound}
\abs{\partial^{\alpha} \delta(x)} \leq C\abs{\delta(x)}^{1-\gamma} \quad \text{for all $x \in \R^n$}.
\end{align}

With this horizon at hand, we define a heterogeneous kernel function
\begin{equation}\label{eq:rhodef}
\rho:\Omega \times \R^n\setminus \{0\} \to [0,\infty), \quad \rho(x,z):=\d(x)^{-n}\rho_1\left(\frac{z}{\d(x)}\right).
\end{equation}
Notice that $\rho(x,\cdot)$ as a rescaled version of $\rho_1$ by a factor of $\d(x)$ satisfies
\begin{equation}\label{eq:rhonormalization}
\supp(\rho(x,\cdot)) \subset B_{\d(x)}(0) \quad \text{and} \quad \norm{\rho(x,\cdot)}_{L^{1}(\R^n)}=\norm{\rho_1}_{L^1(\R^n)}=n \quad \text{for all $x \in \Omega$}.
\end{equation}

The nonlocal gradient with heterogeneous horizon for smooth functions is now defined as follows. Later, we will be able to extend the definition to integrable functions and, more broadly, to distributions (cf.~\eqref{DrviaQro}). 
\begin{definition}[Heterogeneous nonlocal gradient I]\label{def:inhomgradient}
For $\phi \in \Scal(\R^n)$,
we define the heterogeneous nonlocal gradient $\Gr \phi:\Omega \to \R^n$ as
\begin{align}\label{def:nonlocalgradient_sec2}
\Gr \phi(x) = \int_{\Omega} \frac{\phi(y)-\phi(x)}{\abs{y-x}}\frac{y-x}{\abs{y-x}}\rho(x,y-x)\,dy \qquad\text{ for $x \in \Omega$,}
\end{align}
with $\rho$ the kernel function from~\eqref{eq:rhodef}. 
\end{definition}

Observe that the integral in~\eqref{def:nonlocalgradient_sec2} is convergent for all $x \in \Omega$ due to \eqref{eq:rhonormalization} and gives rise to a bounded function on $\Omega$.

We point out that $\Gr$ mirrors the structure of the homogeneous nonlocal gradient in Definition~\ref{def:nonlocalgradient}, but has a kernel function depending also on the $x$-variable through \eqref{eq:rhodef}; our notation, which makes a slight abuse of the fact that $\rho$ depends on two variables, is chosen to represent this property. Intuitively speaking, the gradient $\Gr$ at $x \in \Omega$ agrees with the nonlocal gradient $D_{\rho_1}$ after rescaling by $\d(x)$, and as $\dist(x,\partial\Omega)\to 0$, should approximate the classical local gradient due to the vanishing of the horizon (cf.~\cite{MeS,CKS25,Arr25} on localization in the homogeneous setting).

 An alternative characterization of $\Gr$ via the function $Q_{\rho_1}$ from \eqref{eq:Qrhodefinition} is possible as well, which illuminates the connection with the theory of pseudo-differential operators. Indeed, we note that rescaling \eqref{eq:translation} by the position-dependent horizon $\d(x)$ yields for $\phi \in \Scal(\R^n)$ that
\[
\Gr \phi (x) =\int_{ \Omega}\d(x)^{-n}Q_{\rho_1}\left(\frac{x-y}{\d(x)}\right)\nabla \phi(y)\,dy \quad \text{for all $x \in \Omega$},
\]
or expressed via the Fourier transform,
\begin{equation}\label{eq:pdoform}
\Gr \phi (x) = \int_{\R^n} e^{2\pi i x \cdot\xi} \widehat{Q}_{\rho_1}(\d(x)\xi) \, 2\pi i \xi\,\widehat{\phi}(\xi)\,d\xi \quad \text{for all $x \in \Omega$};
\end{equation}
here, we have used that the Fourier transform of $z\mapsto \d(x)^{-n}Q_{\rho_1}(z/\d(x))$ is equal to $\xi \mapsto\widehat{Q}_{\rho_1}(\d(x)\xi)$, given the scaling properties of the Fourier transform.  From \eqref{eq:pdoform}, one can see that $D_{\rho(\cdot)}$ has the form of a pseudo-differential operator, at least on the set $\Omega$. In fact, by   defining the globally smooth symbol
\begin{align}\label{symbol_qrho}
q_{\rho(\cdot)}:\R^n \times \R^n \to \R, \quad q_{\rho(\cdot)}(x,\xi):=\widehat{Q}_{\rho_1}(\d(x)\xi) \quad 
\end{align}
and its associated pseudo-differential operator
\begin{align*}
 \Qcal_{\rho(\cdot)}:=\Op(q_{\rho(\cdot)}),
 \end{align*}
 \color{black}
we find in view of $q_{\rho(\cdot)}=\widehat{Q}_{\rho_1}(0) = \norm{Q_{\rho_1}}_{L^1(\R^n)}=1$ on $\Omega^c\times \R^n$ that for $\phi \in \Scal(\R^n)$,
\begin{equation}\label{eq:Qrrepr}
\Qr \phi (x) = \begin{cases}
\displaystyle \int_{ \Omega}\d(x)^{-n}Q_{\rho_1}\left(\frac{x-y}{\d(x)}\right) \phi(y)\,dy &\text{if $x \in \Omega$},\\
 \phi(x) \quad &\text{if $x \in \Omega^c$},
\end{cases} 
\end{equation}
and thus, also
\begin{equation}\label{eq:Drcoupling}
\Qr \nabla \phi (x) = \begin{cases}
\Gr \phi (x) \quad &\text{if $x \in \Omega$},\\
\nabla \phi(x) \quad &\text{if $x \in \Omega^c$}.
\end{cases}
\end{equation}
 In this way, the operator $\Qr$ induces a nonlocal-to-local coupling between the heterogeneous nonlocal and the classical gradient.  

Let us investigate the properties of the operator $\Qr$ more closely. The next proposition identifies $\qr$ is an element of the H\"{o}rmander symbol class, which opens up the possibility of exploiting the well-established results on pseudo-differential operators such as mapping properties.
 
 \begin{proposition}[$\Qr$ as pseudo-differential operator in the H\"ormander class]\label{prop:symbolbounds}
For any $\mu \in (0,1)$, it holds that $\qr \in S^{0}_{1,\mu}$, or equivalently, 
 \begin{align*}
 \Qr=\Op(\qr) \in \OPS^{0}_{1,\mu}.
 \end{align*} 
\end{proposition}

\begin{proof}
The proof follows immediately from  the first part of Lemma~\ref{le:symbolbounds}, considering that $\qr$ is bounded due to the boundedness of $\widehat{Q}_{\rho_1}$.
\end{proof}

While the following lemma is substantial for the previous proposition, it also gives further insights. Specifically, it shows that $\qr$ satisfies additional conditions that will allow us conclude the existence of a parametrix for $\Qr$, see~Section~\ref{sec:operators}\, b). 
\begin{lemma}
\label{le:symbolbounds}
For any $\mu\in (0,1)$ and every $\alpha,\beta \in \N_0^{n}$ there is a constant $C_{\alpha, \beta}>0$ such that
\begin{align}\label{est85}
\abs{\partial^\beta_x \partial^{\alpha}_\xi \qr(x,\xi)} \leq C_{\alpha,\beta} \abss{\xi}^{-\abs{\alpha}+\mu\abs{\beta}}\abs{\qr(x,\xi)} \quad \text{for all $x,\xi\in \R^n$}.
\end{align}
In addition, there exists a constant $C>0$ such that
\begin{align}\label{est86}
\qr(x, \xi) \geq C\langle\xi \rangle^{\lambda-1} \quad\text{ for all $x, \xi\in \R^n$, }
\end{align}
where $\lambda$ is the constant from~\ref{itm:lower}.
\end{lemma}

\begin{proof}
Let $\mu\in (0,1)$ and $\alpha, \beta\in \N_0^n$ be fixed. To show~\eqref{est85}, we observe first  that it suffices to prove the stated bounds for $x\in \Omega$ and $\xi\in \R^n$ with $\abs{\xi} \geq 1$, 
due to the smoothness and positivity of $\qr$ along with the fact that $\qr$ is constant on $\Omega^c\times \R^n$.

We can compute that
\[
\partial^{\alpha}_\xi \qr(x,\xi) = \d(x)^{\abs{\alpha}}\partial^{\alpha}\widehat{Q}_{\rho_1}(\d(x)\xi).
\]
 Using the bounds in~\eqref{est2:lemmaBeMCSc24} of Proposition~\ref{prop:BeMCSc24} and \eqref{eq:abssproperty} already gives us the desired estimates when $\beta=0$.

 For the case $\beta \not =0$, we find with the Leibniz rule that 
\begin{align*}
\partial^\beta_x \partial^{\alpha}_\xi \qr(x,\xi) 
=  \sum_{\beta'\in \N_0^n, \,\beta'\leq \beta} {\frac{\beta!}{\beta'!(\beta-\beta')!}} \bigl(\partial_x^{\beta-\beta'} \delta(x)^{|\alpha|}\bigr) \bigl( \partial_x^{\beta'}\partial^\alpha \widehat{Q}_{\rho_1}(\d(x)\xi)\bigr).
\end{align*}
The chain rule along with the decay behavior of $\delta$ from~\eqref{eq:decaybound}  with some $\gamma\in (0,1)$ (to be determined below) allows us to deduce for any $\beta'\in \N_0^n$ that
\begin{align*}
\abs{ \partial_x^{\beta-\beta'} \delta(x)^{|\alpha|} } \leq C \sum_{j=0}^{\min\{|\beta-\beta'|, |\alpha|\}} |\delta(x)|^{|\alpha|-j} |\delta(x)|^{j(1-\gamma)} \leq C  |\delta(x)|^{|\alpha|(1-\gamma)}
\end{align*}
and
\begin{align*}
 \abs{\partial_x^{\beta'} \partial^\alpha \widehat{Q}_{\rho_1}(\d(x)\xi)}  \leq C \sum_{k=0}^{|\beta'|} |\delta(x)|^{k(1-\gamma)} |\xi|^k \abs{\nabla^k \partial^\alpha \widehat{Q}_{\rho_1}(\d(x)\xi)}
\end{align*}
with constants $C>0$. Altogether,  it follows then in view of~\eqref{est2:lemmaBeMCSc24} of Proposition~\ref{prop:BeMCSc24} that 
\begin{align}\label{eq:leibnizterms}
\abs{\partial^\beta_x \partial^{\alpha}_\xi \qr(x,\xi)} 
 \leq C \sum_{k=0}^{\abs{\beta}} A_k(x, \xi)\abs{\qr(x,\xi)},
\end{align}
with $A_{k}(x, \xi) :=  \abs{\d(x)}^{(\abs{\alpha}+k)(1-\gamma)}\abs{\xi}^{k} \abss{\d(x)\xi}^{-\abs{\alpha}-k}$. 

As we will see below, choosing $\gamma$ such that $\gamma(\abs{\alpha}+\abs{\beta})\leq \mu \abs{\beta}$, which is possible since  $\abs{\beta} \geq 1$, finally yields~\eqref{est85}. Indeed, for any $k\in \{0, \ldots, |\beta|\}$, we can estimate in the case $\d(x) \geq \abs{\xi}^{-1}$ that 
\begin{align*}
A_{k}(x, \xi) &\leq \abs{\d(x)}^{(\abs{\alpha}+k)(1-\gamma)}\abs{\xi}^{k} \abs{\d(x)\xi}^{-\abs{a}-k}\\
&\leq \abs{\d(x)}^{-\gamma(\abs{\alpha}+k)} \abs{\xi}^{-\abs{\alpha}}\leq \abs{\xi}^{-\abs{\alpha}+\gamma(\abs{\alpha}+k)} \leq C\abss{\xi}^{-\abs{\alpha}+\mu\abs{\beta}},
\end{align*}
and if $\d(x) \leq \abs{\xi}^{-1}$, one has
\begin{align*}
A_{k}(x, \xi) &\leq \abs{\d(x)}^{(\abs{\alpha}+k)(1-\gamma)}\abs{\xi}^{k} \leq \abs{\xi}^{-(\abs{\alpha}+k)(1-\gamma)+k}=\abs{\xi}^{-\abs{\alpha}+\gamma(\abs{\alpha}+k)} \leq C\abss{\xi}^{-\abs{\alpha}+\mu\abs{\beta}}.
\end{align*}
Filling these bounds into \eqref{eq:leibnizterms} gives us the desired estimates.

 To prove~\eqref{est86}, we can use Proposition~\ref{prop:BeMCSc24} and \ref{itm:lower} to find that
\begin{align*}
\widehat{Q}_{\rho_1}(\zeta) \geq C \abss{\zeta}^{1-\lambda} \quad \text{for all $\zeta \in \R^n$.}
\end{align*}
Hence, we find that
\[
\qr(x, \xi) =\widehat{Q}_{\rho_1}(\delta(x)\xi)\geq C \abss{\delta(x)\xi}^{\lambda-1}\geq C\abss{\xi}^{\lambda-1} \quad \text{for all $x,\xi \in \R^n$.}\qedhere
\] 
\black
\end{proof}

As a consequence of~Proposition~\ref{prop:symbolbounds}, $\Qr$ maps $\Scal(\R^n)$ continuously into itself and extends by duality to $\Scal'(\R^n)$. This extension lies also in $\LBH^0$, and is thus, in particular, bounded on $L^p(\R^n)$ for all $p \in (1,\infty)$, cf.~Section~\ref{sec:pdo}. Moreover, we have due to~\eqref{eq:Qrrepr}, the density of smooth test functions in $\Scal'(\R^n)$ and the continuity of the restriction operator $\pi_{\overline{\Omega}^c}:\Dcal'(\R^n) \to \Dcal'(\overline{\Omega}^c)$ in \eqref{eq:piOmega} that
\begin{align}\label{Qridentityoutside}
{(\Qr u ) |_{\overline{\Omega}^c} = u |_{\overline{\Omega}^c}} \quad \text{for all $u \in \Scal'(\R^n)$.}
\end{align}

We will establish next that $\Qr$ admits a restriction to $\Omega$ in the sense of Section~\ref{sec:pdo}\,c). 
To see this, it suffices to show that the dual operator $\Qr^\ast$, which belongs to $\OPS^{0}_{1,\mu}$ for all $\mu >0$, maps $C_c^{\infty}(\Omega)$ continuously to itself. Indeed,  in view of~\eqref{eq:Qrrepr}, a simple application of Fubini's theorem for $\phi, \psi \in C_c^{\infty}(\R^n)$ implies 
\begin{align*}
\int_{\R^n} \Qr \phi \,\psi\,dx&=\int_{\Omega}\int_{\Omega} \d(x)^{-n} Q_{\rho_1}\left(\frac{x-y}{\d(x)}\right)\phi(y)\,dy \,\psi(x)\,dx + \int_{\Omega^c} \phi\,\psi\,dx \\
&= \int_{\Omega} \int_{\Omega} \d(y)^{-n} Q_{\rho_1}\left(\frac{x-y}{\d(y)}\right)\psi(y)\,dy \,\varphi(x)\,dx+ \int_{\Omega^c} \psi\,\phi\,dx,
\end{align*}
which yields the explicit formula
\begin{align}\label{dualQast}
\Qr^* \psi(x) = \begin{cases}
\displaystyle \int_{\Omega} \d(y)^{-n} Q_{\rho_1}\left(\frac{x-y}{\d(y)}\right)\psi(y)\,dy &\text{if $x \in \Omega$},\\
\psi(x) &\text{if $x \in\Omega^c$}.
\end{cases}
\end{align}

This representation of $\Qr^*\in \OPS^{0}_{1,\mu}$ allows us to conclude that $\Qr^*$ is a continuous map from $C_c^{\infty}(\R^n)$ and $C_c^{\infty}(\Omega)$, respectively, into the same spaces, as desired.  Hence, we can define $\Qr$ on the space of distributions $\Dcal'(\Omega)$ as
\begin{align}\label{eq:Qrrestriction}
\Qro:\Dcal'(\Omega) \to \Dcal'(\Omega),\, \quad \inner{\Qro u, \psi} = \inner{u,\Qr^* \psi} \quad \text{for $\psi \in C_c^{\infty}(\Omega)$},
\end{align}
cf.~\eqref{eqrestrict}. From the abstract theory, we know also that 
$\Qro$ is bounded on $H^{s,p}(\Omega)$ for all $s \in \R$, or in other words, 
\begin{align}\label{QroLBH0}
\Qro\in \LBH_\Omega^0.
\end{align} 

With the restriction operator $\pi_\Omega$ from~\eqref{eq:piOmega}, one finds for all $u\in \Scal'(\R^n)$ that
\begin{align}\label{Qcal_restriction}
\Qro  (\pi_\Omega u) = \pi_\Omega(\Qr u). 
\end{align}
Clearly, $\Qro$ maps $C_c^\infty(\Omega)$ to $C_c^\infty(\Omega)$ by \eqref{eq:Qrrepr}, and if $\varphi \in C^\infty(\overline{\Omega})$, we  have that $\Qro \phi\in C^\infty(\overline{\Omega})$ agrees with the integral in the first case in \eqref{eq:Qrrepr}, i.e.,
\begin{align*}
\Qro\varphi(x) = \int_{\Omega}\d(x)^{-n}Q_{\rho_1}\left(\frac{x-y}{\d(x)}\right) \phi(y)\,dy \quad \text{for $x \in \Omega$}. 
\end{align*}

\medskip

 After examining the operators $\Qr$ and $\Qro$, let us shift our focus back to the nonlocal gradient objects we are mainly interested in. 
 The observation from~\eqref{eq:Drcoupling} that $\Gr \varphi =\Qr \nabla \varphi$ on $\Omega$ for any Schwartz function $\varphi$ suggests a  natural way to extend the nonlocal gradient $\Gr$ to the space of distributions, as done in the next definition.  There,  we also introduce a whole-space analogue of $\Gr$, cf.~\eqref{eq:Drcoupling}.

\begin{definition}[Heterogeneous nonlocal gradients II]\label{def:hetnonlocgrad2}
For $u\in\Dcal'(\Omega)$, we define its heterogeneous nonlocal gradient by 
\begin{align}\label{DrviaQro}
\Gr u:=\Qro\nabla u.
\end{align} 
Further, we introduce the whole-space heterogeneous nonlocal gradient on $\Scal'(\R^n)$ as the operator $\Grr:=\Qr \nabla$. 
\end{definition}
Note that in view of the properties of $\Qr$ and $\Qro$, it holds that
$$\Gr\in \LBH^1_\Omega\qquad \text{and }\qquad\Grr \in 
 \LBH^1,$$ 
respectively, cf.~\eqref{prop:intbyparts} and~Proposition~\ref{prop:symbolbounds}. Since
\begin{align*}
\inner{\Gr u,\psi} &= \inner{\Qro \nabla u,\psi}=-\inner{u,\Div \Qr^{*}\psi} \quad \text{for $u\in \D'(\Omega)$ and $\psi \in C_c^{\infty}(\Omega;\R^n)$, }\\
\inner{\Grr u,\psi} &= \inner{\Qr \nabla u,\psi}=-\inner{u,\Div \Qr^{*}\psi} \quad \text{for $u\in \Scal'(\R^n)$ and $\psi \in \Scal(\R^n;\R^n)$,}
\end{align*}
 it is straightforward to infer 
 that $\Gr =\Grr^\Omega$.
Moreover, these identities~
show that $$\Div^*_{\rho(\cdot)} := \Div  \Qr^*$$ acts as a dual operator to $\Gr$ and $\Grr$. In fact, we can prove an integration by parts formula with local boundary conditions where this operator arises naturally; for a generalization see~\eqref{IbP_general}. 
 
\begin{proposition}[Integration by parts formula for smooth functions]\label{prop:intbyparts} 
Let $\phi \in C_c^{\infty}(\R^n)$ and $\psi \in C_c^{\infty}(\R^n;\R^n)$, then it holds that
\[
\int_{\Omega} \Gr \phi \cdot \psi\,dx = -\int_{\Omega} \phi\,\Div^*_{\rho(\cdot)}\psi\,dx+\int_{\partial \Omega} \phi\,\psi\cdot\nu \,d\Hcal^{n-1},
\]
with $\nu$ an outward pointing unit normal to $\partial \Omega$.
\end{proposition}
\begin{proof} With $\phi \in C_c^{\infty}(\R^n)$ and $\psi \in C_c^{\infty}(\R^n;\R^n)$, 
we can compute using \eqref{eq:Drcoupling} that
\begin{align*}
\int_{\Omega} \Gr \phi \cdot \psi\,dx &+\int_{\Omega^c}\nabla \phi \cdot \psi\,dx= \int_{\R^n} \Qr \nabla\phi \cdot \psi\,dx\\
&=\int_{\R^n} \nabla \phi \cdot \Qr^*\psi\,dx= \int_{\Omega} \nabla \phi \cdot \Qr^* \psi\,dx+\int_{\Omega^c} \nabla \phi \cdot \psi\,dx\\
&= -\int_{\Omega} \phi \, \Div \Qr^* \psi\,dx +\int_{\partial \Omega} \phi\,\Qr^*\psi\cdot\nu \,d\Hcal^{n-1} +\int_{\Omega^c} \nabla \phi \cdot \psi\,dx \\
&= -\int_{\Omega} \phi \, \Div^*_{\rho(\cdot)}\psi\,dx +\int_{\partial \Omega} \phi\,\psi\cdot\nu \,d\Hcal^{n-1} +\int_{\Omega^c} \nabla \phi \cdot \psi\,dx;
\end{align*}
in the third line, we have used the classical integration by parts formula, given that $\Qr^*\psi$ is smooth, whereas the last line uses that $\Qr^*\psi=\psi$ on $\partial \Omega$. Finally, subtracting $\int_{\Omega^c} \nabla \phi \cdot \psi\,dx$ from both sides yields the desired result.
\end{proof}

\subsection{Auxiliary  operators related to the heterogeneous nonlocal gradient }\label{sec:operators} 
 
Next, we introduce several important operators closely related to $\Qr$ and $\Gr$, which play a central role throughout the paper and will be extensively utilized in the subsequent analysis. These include commutators, parametrices, and  residuals. 
In the case of bounded domains, analogous operators can be defined either by restricting pseudo-differential operators, as outlined in Section~\ref{sec:pdo}\,c), or through a combination of suitable extensions and restrictions. 
 \smallskip

\textit{a) Commutators.}  To estimate the difference between $\Grr$ and $\nabla \Qr$, we introduce
\begin{align*}
\Cr:=[\nabla,\Qr] = \nabla\Qr-\Qr \nabla = \nabla \Qr-\Grr \in \OPS^{\mu}_{1,\mu},
\end{align*}
whose symbol is given by $\sigma(\Cr)=i\nabla_x \qr$ and lies in $S^{\mu}_{1,\mu}$ for any $\mu \in (0,1)$; indeed, the latter follows from~\cite[Theorem~18.1.6]{Hor07}, considering that  $\qr\in S^0_{1, \mu}$ for any $\mu\in (0,1)$ by Lemma~\ref{le:symbolbounds}.

Using the restricted operator $\Qr^\Omega$, one can define in the same way
\begin{align}\label{eq:cro}
\Ccal_{\rho(\cdot)}^{\Omega} :=  [\nabla,\Qro] = \nabla\Qro-\Qro \nabla=\nabla\Qro-\Gr. 
\end{align}
It is straightforward to show that this coincides with the restriction of $\Cr$ in the sense of~\eqref{mappingproperties_AcalOmega}.

Altogether, we have the commutators
\begin{align}\label{mappingprop_C}
\Cr\in \LBH^+\quad \text{and} \quad  \Cr^\Omega\in \LBH_\Omega^+,\end{align}
in the whole space and bounded setting, respectively.

\smallskip

\textit{b) Parametrix and residuals.} We observe in light of Lemma~\ref{le:symbolbounds}  that 
$\qr$ satisfies~\eqref{eq:parametrixcond2} as well as \eqref{eq:parametrixcond1} with $m'=\lambda-1$; recall that $\lambda\in (0,1)$ is the parameter appearing in~\ref{itm:lower}.  Hence, as outlined in Section~\ref{sec:pdo}\,c), there exists a parametrix symbol $\pr \in S^{1-\lambda}_{1,\mu}$, and since $\qr=1$ on $\Omega^c\times \R^n$, it can be assumed that also 
\begin{equation}\label{eq:p1outside}
\pr(x,\xi)=1 \quad  \text{for all $x \in \Omega^c$ and $\xi\in \R^n$,}
\end{equation} 
given the construction method of the parametrix in \cite[Theorem~5.4]{Kum81}.

The associated operator
\begin{align}\label{eq:pr}
\PR:=\Op(\pr) \in \OPS^{1-\lambda}_{1,\mu}
\end{align}
then satisfies, in particular, $\PR\in \LBH^{1-\lambda}$, meaning
\begin{align}\label{Pbounded}
\PR: H^{s, p}(\R^n)\to H^{s+\lambda-1, p}(\R^n)\quad  \text{is bounded for all $s\in \R$ and $p\in (1, \infty)$.}
\end{align}
Given \eqref{eq:p1outside}, it follows similarly to \eqref{Qridentityoutside} that 
\begin{align}\label{Pridentityoutside}
{(\PR u ) |_{\overline{\Omega}^c} = u |_{\overline{\Omega}^c}} \quad \text{for all $u \in \Scal'(\R^n)$.}
\end{align}

Further, the property of the parametrix yields that the residuals
\begin{align}\label{eq:res}
\Rr:= \PR \Qr - \Id \in \LBH^{-\infty} \quad \text{and} \quad \Rrr:=\Qr \PR-\Id  \in \LBH^{-\infty}
\end{align}
are smoothing operators, which gives a suitable replacement for an exact inverse. Along with \eqref{Qridentityoutside} and \eqref{Pridentityoutside}, we find that
\begin{align}\label{Rprop_complementary}
\text{$\Rr$ and $\Rrr$ map $H^{s, p}(\R^n)$ into  $\Cinfo$ for all $s\in \R$ and $p\in (1, \infty)$.} 
\end{align}
\color{black}
\smallskip

\textit{c) Extension operator.} Let $E$ be the Rychkov universal extension operator, see~\eqref{Rychkov}. 
For $u\in H^{s,p}(\Omega)$ with $s\geq 0$ and $p\in (1, \infty)$, we define
\begin{align}\label{eq:extension}
E_{\rho(\cdot)}u = \begin{cases} u &\text{on $\Omega$},\\
E\Qro u &\text{on $\Omega^c$,}
\end{cases}
\end{align}
which is an element of $L^p(\R^n)$. 
It can be deduced from \eqref{eq:Qrrepr} along with a density argument that for $u\in H^{s,p}(\Omega)$,
\begin{align}\label{QEswitch}
\Qr \Er u = E \Qro u\in H^{s, p}(\R^n). 
\end{align}

We show that $\Er$ is a bounded linear operator $H^{s,p}(\Omega)\to H^{s-1+\lambda, p}(\R^n)$ for any $s\geq 0$ and $p\in (1, \infty)$. Indeed, one can use the mapping properties of $\PR$ from \eqref{eq:pr} along with those of the residuals from \eqref{eq:res} to conclude
\begin{align*}
\norm{\Er u}_{H^{s-1+\lambda,p}(\R^n)}&=\norm{\PR (\Qr \Er u)-\Rr \Er u}_{H^{s-1+\lambda,p}(\R^n)}\\
&\leq C \bigl(\norm{\Qr \Er u}_{H^{s,p}(\R^n)}+C\norm{\Er u}_{L^p(\R^n)}\bigr) \\
&\leq C \bigl(\norm{E(\Qro u)}_{H^{s,p}(\R^n)}+C\norm{u}_{L^p(\Omega)}+C\norm{E(\Qro u)}_{L^p(\Omega^c)}\bigr) \\
&\leq C\bigl(\norm{\Qro u}_{H^{s,p}(\Omega)}+C\norm{u}_{L^p(\Omega)} \bigr)\\
&\leq C\norm{u}_{H^{s,p}(\Omega)};
\end{align*}
more precisely, the second estimate uses that $\PR\in \LBH^{1-\lambda}$ and that $\Rr$ is smoothing and thus, bounded from $L^p(\R^n)$ to $H^{s-1+\lambda,p}(\R^n)$, while the last inequality exploits  the boundedness of $\Qro$ on $H^{s, p}(\Omega)$, cf.~\eqref{QroLBH0}. \smallskip


\textit{d) Restriction of parametrix and residuals.}
We define the bounded linear operator 
\begin{align*}\label{defPrho_parametrix}
\PROmega:= \pi_\Omega \PR E, 
\end{align*} 
where $E$ is again the Rychkov extension operator and $\pi_\Omega$ is the  restriction operator for distributions, cf.~\eqref{eq:piOmega}. It immediately follows from the mapping properties of the operators involved and \eqref{Pridentityoutside} that
\begin{equation}\label{POmegasmooth}
\text{$\PROmega$ maps $\Cinf$ into itself and $C_c^{\infty}(\Omega)$ into $\Cinfo$.}
\end{equation}
and that
\begin{align}\label{POmega}
\PROmega   \in \LBH_\Omega^{1-\lambda};
\end{align}
see~\eqref{Rychkov} and~\eqref{Pbounded}, and note that $\pi_\Omega$ is a bounded linear operator $H^{s, p}(\R^n)\to H^{s, p}(\Omega)$ for all $s\in \R$ and $p\in (1, \infty)$. 

 Together with~\eqref{QEswitch} and~\eqref{eq:res}, we can then deduce for $u\in H^{s,p}(\Omega)$ with $s\geq 0$ and $p\in (1, \infty)$ that
\begin{align*}
\PROmega\Qr^\Omega u & = \pi_\Omega \PR  E \Qr^\Omega  u = \pi_\Omega  \PR \Qr  \Er u \\ & = \pi_\Omega   ({\Id+\Rr)   \Er  u= u+\RrOmega u},
\end{align*}
with $\RrOmega:=\pi_\Omega \Rr \Er$. In view of~the properties of $\Er$ (see Section~\ref{sec:operators}\,c)) and the smoothing behavior of $\Rr$ (see~\eqref{eq:res}), it holds that $\RrOmega$ maps continuously from $H^{s,p}(\Omega)$ to $H^{s+\gamma, p}(\Omega)$ for all $s\geq 0$, $\gamma\in \R$ and $p\in (1, \infty)$, or in short, $ \RrOmega\in \LBH^{-\infty}_{\Omega, 0}$. 

On the other hand, it holds for $u\in H^{s, p}(\Omega)$ with $s\in \R$ and $p\in (1, \infty)$ that
\begin{align*}
\Qr^\Omega  \PROmega u &= \Qr^\Omega   \pi_\Omega    \PR   E u = \pi_\Omega  \Qr  \PR   E u \\ &= \pi_\Omega   {(\Id+\Rrr)   E u=u+ \RrrOmega u,}
\end{align*} 
where  the second identity is due to $\Qr^\Omega   \pi_\Omega \tilde u  = \pi_\Omega  \Qr \tilde u$ for any $\tilde u\in H^{s,p}(\R^n)$ according to~\eqref{Qcal_restriction}, and $\RrrOmega:=\pi_\Omega   \Rrr   E\in \LBH_\Omega^{-\infty}$, cf.~\eqref{eq:res} and~\eqref{Rychkov}. 

In summary, the above demonstrates that the parametrix $\PROmega\in \LBH_\Omega^{1-\lambda}$ and the residual operators
\begin{align}\label{residual_omega}
\RrOmega =\PROmega   \Qr^\Omega -\Id  \in   { \LBH^{-\infty}_{\Omega, 0}} \quad \text{and} \quad \RrrOmega=\Qr^\Omega  \PROmega - \Id \in\LBH_\Omega^{-\infty} 
\end{align}
constitute the counterparts to $\PR$, $\Rr$ and $\Rrr$ from  part b) in the bounded setting.

\begin{remark}
Note that the operator $\PROmega$ is not defined via restriction of the pseudo-differential operator $\PR$ in the sense of Section~\ref{sec:pdo}\,c). Indeed, since it is not clear whether the dual $\PR^\ast$ maps $C_c^\infty(\Omega)$ into itself, the existence of $\PR^\Omega$ is not easily seen.
\end{remark} 

\section{Heterogeneous nonlocal Sobolev spaces}\label{sec:heterogeneous_Sobolevspaces}
We are now in the position to introduce new Sobolev spaces associated to the nonlocal gradient with position-dependent horizon and prove some basic properties. Our results include embedding, density and regularity results for these spaces, as well as the existence of trace and extension operators, and Poincar{\'e} inequalities.  A useful tool is a translation mechanism that allows to switch between classical and nonlocal heterogeneous gradients up to an controlled error.

Let $p \in (1,\infty)$. Recalling the definition of $\Gr$ (see Definitions~\ref{def:inhomgradient} and~\ref{def:hetnonlocgrad2}) and that $\rho_1$ satisfies all the properties \ref{itm:h0}-\ref{itm:upper}, we define the heterogeneous nonlocal Sobolev spaces for functions on bounded, open sets as follows. A version for functions defined on all of $\R^n$ can be found in Section~\ref{subsec:extensions_heterogeneous}. 
\color{black}
\begin{definition}[Heterogeneous nonlocal Sobolev spaces]\label{def:newspaces}
With the nonlocal gradient operator $\Gr=\Qro\nabla$, we define
\[
H^{\rho(\cdot),p}(\Omega):=\{u \in L^p(\Omega) \,:\, \Gr u \in L^p(\Omega;\R^n)\},
\]
endowed with the norm 
\[
\norm{u}_{\Hr}:=\norm{u}_{L^p(\Omega)} +\norm{\Gr u}_{L^p(\Omega;\R^n)}
\]
for $u\in \Hr$.
\end{definition}

Standard techniques show that $\Hr$ are reflexive, separable Banach spaces. As we establish next, the classical Sobolev spaces embed continuously into $\Hr$, which, in turn, are continuously embedded in Bessel potential spaces. 
Recall that $\lambda\in (0,1)$ is the parameter appearing in \ref{itm:lower}.

\begin{proposition}[Embeddings involving $\Hr$]\label{prop:clasembedding}
It holds that 
\begin{align*}
W^{1,p}(\Omega) \hookrightarrow \Hr \hookrightarrow H^{\lambda,p}(\Omega).
\end{align*}  
\end{proposition}

\begin{proof} Let $u\in W^{1,p}(\Omega)$. Due to the boundedness of $\Qro:L^p(\Omega) \to L^p(\Omega)$, cf.~\eqref{QroLBH0}, we find that $\Gr u = \Qro \nabla u \in L^p(\Omega;\R^n)$ with
\[
\norm{\Gr u}_{L^p(\Omega;\R^n)}=\norm{\Qro \nabla u}_{L^p(\Omega;\R^n)} \leq C \norm{\nabla u}_{L^p(\Omega;\R^n)}, 
\]
where $C>0$ is independent of $u$. This yields the first embedding.

To establish the second one, we first prove the weaker statement that
\begin{align}\label{embedding_firststep}
\Hr\hookrightarrow H^{\lambda-\mu, p}(\Omega)\quad\text{ for all $\mu >0$.}
\end{align}
Indeed, by exploiting the mapping properties of the operators from Section~\ref{sec:operators}\, a)  and~c) as well as the equivalence of norms in~\eqref{eq:besselboundednorm}, we find that for all $u \in \Hr$,
\begin{align}\label{est67}
\begin{split}
\norm{ u}_{H^{\lambda-\mu,p}(\Omega)} &=\norm{\PROmega (\Qr^\Omega u) -\RrOmega u}_{H^{\lambda-\mu,p}(\Omega)}\\
&\leq C\bigl(\norm{\Qr^\Omega u}_{H^{1-\mu, p}(\Omega)} + \norm{u}_{L^p(\Omega)}\bigr)\\
& \leq C\bigl(\norm{\Qr^\Omega u}_{H^{-\mu, p}(\Omega)} + \norm{\nabla \Qr^\Omega u}_{H^{-\mu, p}(\Omega;\R^n)} + \norm{u}_{L^p(\Omega)} \bigr)\\
& \leq C\bigl(\norm{u}_{H^{-\mu, p}(\Omega)} + \norm{\Cr^{\Omega} u}_{H^{-\mu, p}(\Omega;\R^n)} + \norm{\Gr u}_{H^{-\mu, p}(\Omega;\R^n)}+  \norm{u}_{L^p(\Omega)} \bigr)\\
& \leq C\bigl( \norm{u}_{\Hr} + \norm{\Cr^{\Omega} u}_{H^{-\mu, p}(\Omega;\R^n)} \bigr)
\end{split}
\end{align}
with a generic constant $C>0$. 
As $\Cr^{\Omega}: L^p(\Omega)\to H^{-\mu, p}(\Omega;\R^n)$ is bounded in view of~\eqref{mappingprop_C}, we obtain further that
\begin{align*}
\norm{u}_{H^{\lambda-\mu,p}(\Omega)} \leq C \norm{u}_{\Hr} \quad \text{for all $u \in \Hr$, }
\end{align*}
which implies~\eqref{embedding_firststep}.

 Now we repeat the same estimate of~\eqref{est67} to show that 
\begin{align*}
\norm{u}_{H^{\lambda,p}(\Omega)} \leq  C\bigl( \norm{u}_{\Hr} + \norm{\Cr^{\Omega} u}_{L^p(\Omega;\R^n)} \bigr) \quad \text{for all $u \in \Hr$.}
\end{align*}
Using that $\Cr^{\Omega}: H^{\lambda-\mu}(\Omega)\to L^p(\Omega;\R^n)$ is bounded for any $\mu\in (0, \lambda)$  along with~\eqref{embedding_firststep} allows us then to conclude $\Hr \hookrightarrow H^{\lambda,p}(\Omega)$, as stated.
\end{proof}

We now present a major tool of this paper, which establishes that $\Qro$ is a bounded linear operator from $\Hr$ into $W^{1,p}(\Omega)$ and that $\PROmega$ acts in the reverse direction. The heterogeneous nonlocal gradient $\Gr$ can be translated into a classical gradient, and vice versa, up to operators of lower order. 
This bidirectional translation plays a crucial role in connecting the nonlocal and classical Sobolev theories.

\begin{lemma}[Translation between heterogeneous nonlocal and classical Sobolev spaces]\label{le:connection}
The following statements hold:
\begin{itemize}
\item[$(i)$] The operator $\Qr^\Omega:\Hr \to W^{1,p}(\Omega)$ is bounded and
\begin{align}\label{translation1}
\nabla \Qr^\Omega -\Gr \in \LBH_\Omega^+. 
\end{align}
\item[$(ii)$] The operator $\PROmega:W^{1,p}(\Omega) \to \Hr$ is bounded and
\begin{align}\label{translation2}
\Gr \PROmega  -\nabla   \in  \LBH_\Omega^{\gamma} \qquad\text{for all $\gamma>1-\lambda$.}
\end{align}
\end{itemize}
\end{lemma}
\begin{proof}  To prove $(i)$, we observe first that~\eqref{translation1} is a direct consequence of the mapping properties of the restricted commutator~$\Cr^{\Omega}=\nabla \Qr^\Omega -\Gr\in \LBH_\Omega^+$, see~\eqref{eq:cro} and~\eqref{mappingprop_C}.

The stated boundedness of $\Qr^\Omega$ as an operator from $\Hr$ to $W^{1,p}(\Omega)$ follows from the fact that $\Qr^\Omega:L^p(\Omega)\to L^p(\Omega)$ is bounded, cf.~\eqref{QroLBH0}, together with the following estimate: For all $u\in \Hr$,
we deduce that
 \begin{align*}
\norm{\nabla \Qro u}_{L^p(\Omega;\R^n)} &\leq \norm{\Gr u}_{L^p(\Omega;\R^n)} + \norm{\Cr^{\Omega}u}_{L^p(\Omega;\R^n)} \\
&\leq \norm{\Gr u}_{L^p(\Omega;\R^n)} + C\norm{u}_{H^{\lambda,p}(\Omega)}\leq C\norm{u}_{\Hr},
\end{align*}
 with a constant $C>0$, where we have employed~\eqref{mappingprop_C} as well as the embedding $\Hr \hookrightarrow H^{\lambda,p}(\Omega)$ from Proposition~\ref{prop:clasembedding}.

For part $(ii)$, we can compute with the restricted operators from Section~\ref{sec:operators}\, a)  and d) that
\begin{align*}
\Gr \PROmega - \nabla &= \Qr^\Omega \nabla \PROmega - \nabla \\
&= \nabla \Qr^\Omega\PROmega -\nabla { - } \Cr^\Omega\PROmega\\ & = \nabla \RrrOmega { - } \Cr^\Omega\PROmega. 
\end{align*} 
Considering that $\Cr^\Omega\in \LBH_\Omega^+$, 
$\PROmega\in \LBH_\Omega^{1-\lambda}$, 
and the smoothing behavior of $\RrrOmega \in \LBH_\Omega^{-\infty}$ by~\eqref{residual_omega} then implies~\eqref{translation2}. 

 Further, we infer for any $u\in W^{1,p}(\Omega)$ that
 \begin{align*}
\norm{\Gr \PROmega u}_{L^p(\Omega;\R^n)} &\leq \norm{\nabla u}_{L^p(\Omega;\R^n)} + \norm{\nabla u - \Gr\PROmega u}_{L^p(\Omega;\R^n)}    \\
&\leq \norm{\nabla u}_{L^p(\Omega;\R^n)} + C\norm{u}_{H^{\gamma, p}(\Omega)}\\ & \leq C\norm{u}_{W^{1,p}(\Omega)},
\end{align*}
with a constant $C>0$; while the second inequality uses~\eqref{translation2} with some $\gamma\in (1-\lambda, 1)$, the third one follows from $W^{1,p}(\Omega)\hookrightarrow H^{\gamma, p}(\Omega)$, cf.~\eqref{Hcompactembedding}. 
Along with the boundedness of $\PROmega\in\LBH_\Omega^{1-\lambda}$ (see~\eqref{POmega}) as linear operator from $W^{1,p}(\Omega)$ to $H^{\lambda, p}(\Omega)\hookrightarrow L^p(\Omega)$, this completes the proof of $(ii)$.
\end{proof}

A simple, yet interesting, observation arising from the translation result in Lemma~\ref{le:connection} is that a function $u \in L^p(\Omega)$ lies in $\Hr$ if and only if $\Qro u$ belongs to the classical Sobolev space $W^{1,p}(\Omega)$. In other words, one can alternatively characterize the heterogeneous nonlocal spaces from Definition~\ref{def:newspaces} as follows.
\begin{corollary}\label{cor:represHrhop}
It holds that
\begin{align*}
\Hr=\{u\in L^p(\Omega): \Qr^\Omega u\in W^{1,p}(\Omega)\},
\end{align*}
and the norm $\norm{\cdot}_{\Hr}$ is equivalent to $\norm{\Qr^\Omega\,\cdot\,}_{W^{1,p}(\Omega)}$. 
\end{corollary}

\subsection{Traces}\label{subsec:traces}
As a further consequence of Lemma~\ref{le:connection}, we will now deduce the existence of a trace operator on $\Hr$, which is one of the main results in this work. In fact, this operator naturally extends the standard trace operator on the Sobolev space $W^{1,p}(\Omega)$, denoted by $T_{W^{1,p}(\Omega)}$ in the following, to our heterogeneous nonlocal spaces. Note that the resulting trace space for $\Hr$ is exactly equal to the classical one of $W^{1,p}(\Omega)$, that is, $W^{1-1/p,p}(\partial \Omega)$. 
\begin{theorem}[Existence of a trace operator]\label{thm:trace}
There exists a unique bounded linear operator $\Trho:\Hr \to W^{1-1/p,p}(\partial \Omega)$ such that
\begin{align}\label{tracecondition}
\Trho u = T_{W^{1,p}(\Omega)} u \qquad \text{for all $u \in W^{1,p}(\Omega)$}.
\end{align}
Moreover, there is a bounded linear operator $\Srho:W^{1-1/p,p}(\partial \Omega) \to \Hr$ with
\begin{align*}
\Trho   \Srho = \Id.
\end{align*}
\end{theorem}

\begin{proof}
The uniqueness follows from the density of $W^{1,p}(\Omega)$ in $\Hr$ (see~Theorem~\ref{th:density}\,$(i)$ below). For the existence, we define $\Trho$ as the composition of the standard trace operator of $W^{1,p}(\Omega)$ with $\Qro$ from~\eqref{eq:Qrrestriction}, that is, 
\begin{align}\label{def:traceoperator}
\Trho u :=T_{W^{1,p}(\Omega)} \Qro u \qquad \text{ for $u \in \Hr$. }
\end{align} The stated mapping properties and the boundedness of $\Trho$ follow immediately from those of $T_{W^{1,p}(\Omega)}$ together with the boundedness of $\Qro:\Hr \to W^{1,p}(\Omega)$ from Lemma~\ref{le:connection}\,$(i)$.

To prove~\eqref{tracecondition}, we start by considering $\tilde \phi \in C_c^{\infty}(\R^n)$ and $\phi=\pi_\Omega \tilde \phi \in C^{\infty}(\overline{\Omega})$ and recall that~\eqref{eq:Qrrepr} implies 
 \begin{align*}
 \Qr  \tilde \phi =\tilde \phi \quad \text{ on $\partial \Omega$.}
 \end{align*} 
 It follows then along with~\eqref{Qcal_restriction} and the properties of the classical trace operator that
\begin{align*}
\Trho \phi =T_{W^{1,p}(\Omega)}\Qro \phi = T_{W^{1,p}(\Omega)} \pi_\Omega  \Qr  \tilde{\phi}  =\tilde \phi|_{\partial \Omega}=T_{W^{1,p}(\Omega)}\phi.
\end{align*}
Now let $u \in W^{1,p}(\Omega)$ and $(\phi_j)_j \subset C^{\infty}(\overline{\Omega})$ converge to $u$ in $W^{1,p}(\Omega)$ as $j\to \infty$. Then, we find
\begin{align*}
T_{W^{1,p}(\Omega)} u=\lim_{j \to \infty} T_{W^{1,p}(\Omega)}\phi_j = \lim_{j \to \infty} \Trho \phi_j=\Trho u,
\end{align*}
where the last convergence uses that $\phi_j \to u$ in $\Hr$ according to Proposition~\ref{prop:clasembedding}. This density argument shows that $\Trho$ in fact extends the trace operator of $W^{1,p}(\Omega)$. 

Lastly, for $\Srho$, one can simply take any linear and continuous right inverse of $T_{W^{1,p}(\Omega)}$, given the continuous embedding $W^{1,p}(\Omega) \hookrightarrow \Hr$ of Proposition~\ref{prop:clasembedding}.
\end{proof}

In light of the previous result, one can define affine subspaces of $\Hr$, which collect functions with prescribed boundary values. 
\begin{definition}[Nonlocal spaces with prescribed trace]\label{def:trace}
For $g \in W^{1-1/p,p}(\partial \Omega)$, let
\[
\Hrg :=\{u \in \Hr \,:\, \Trho (u)=g \ \text{on $\partial \Omega$}\}.
\]
\end{definition}
Notice that these sets are non-empty, given the existence of the right inverse of the trace operator $\Trho$. 
Moreover, the translation procedure preserves traces, i.e., it holds that 
\begin{align*}
\Qro:H_g^{\rho(\cdot), p}(\Omega) \to  W_g^{1,p}(\Omega) \quad \text{and} \quad \PROmega:W_g^{1,p}(\Omega) \to H_g^{\rho(\cdot), p}(\Omega) \ \text{are bounded.}
\end{align*}
Indeed, the first of these observations follows by definition, while for the second we compute for $v \in W^{1,p}(\Omega)$ that
\begin{align*}
\Trho \PROmega v &= T_{W^{1,p}(\Omega)} \Qro \PROmega v \\
&=T_{W^{1,p}(\Omega)}\bigl(v+\RrrOmega v\bigr) = T_{W^{1,p}(\Omega)}v,
\end{align*}
using that $ T_{W^{1,p}(\Omega)}\RrrOmega v =0$ by \eqref{Rprop_complementary} and the definition of $\RrrOmega$. \color{black}

\subsection{Density results}
Let us continue by using the translation mechanism to prove the following useful density results for $\Hr$ and $\Hro$, in the spirit of the Meyers-Serrin theorem for classical Sobolev spaces.

\begin{theorem}[Density of smooth functions]\label{th:density}
These two statements hold:
\begin{itemize}
\item[(i)] The space $\Cinf$ is dense in $\Hr$. \smallskip

\item[(ii)] The space $C_c^{\infty}(\Omega)$ is dense in $\Hro$.
\end{itemize}
\end{theorem}

\begin{proof}
Both parts follow in a similar way by reducing the approximation step to the case of classical Sobolev spaces via the  translation mechanism in
Lemma~\ref{le:connection}.

\smallskip

\textit{Part (i).} Take $u \in \Hr$ and define $v=\Qro u \in W^{1,p}(\Omega)$, considering Lemma~\ref{le:connection}\,$(i)$. Then, there exists a sequence $(\phi_j)_j \subset \Cinf$ with $\phi_j \to v$ in $W^{1,p}(\Omega)$ as $j\to \infty$. We obtain in view of~\eqref{POmegasmooth} that $\psi_j:=\PROmega\phi_j \in \Cinf$  for  all $j \in \N$, and deduce using Lemma~\ref{le:connection}\,$(ii)$
that
\[
\psi_j \to \PROmega v =  u +\RrOmega u \quad \text{in $\Hr$ as $j \to \infty$.}
\]
Since $\RrOmega \in \LBH^{-\infty}_{\Omega, 0}$, it follows that $\RrOmega u \in \Cinf$, so that $(\psi_j-\RrOmega u)_j$  defines a sequence in $\Cinf$ converging to $u$ in $\Hr$. 

\smallskip

\textit{Part (ii).} It suffices to prove that $\Cinfo$ is dense in $\Hro$, given the density of $C_c^{\infty}(\Omega)$ in $\Cinfo$ with respect to the $W^{1,p}(\Omega)$-norm along with the first embedding in Proposition~\ref{prop:clasembedding}. This can be done in exactly the same manner as part $(i)$, by taking $(\phi_j)_j \subset C_c^{\infty}(\Omega)$ and using the mappings properties $\PROmega:C_c^{\infty}(\Omega) \to \Cinfo$ by \eqref{POmegasmooth} and $\RrOmega:L^p(\Omega) \to \Cinfo$ from \eqref{Rprop_complementary}.
\end{proof}
To close this subsection, we remark that combining the density result of Theorem~\ref{th:density} with Theorem~\ref{thm:trace} allows us to extend the integration by parts formula from Proposition~\ref{prop:intbyparts} to $\Hr$.

\begin{corollary}[Integration by parts formula] For all $u \in \Hr$ and $\psi \in C_c^{\infty}(\R^n;\R^n)$, it holds that
\begin{align}\label{IbP_general}
\int_{\Omega} \Gr u \cdot \psi\,dx = -\int_{\Omega} u\,\Div^*_{\rho(\cdot)}\psi\,dx+\int_{\partial \Omega} \Trho u \ \psi\cdot\nu \,d\Hcal^{n-1}. 
\end{align}
\end{corollary}

\subsection{Extensions}\label{subsec:extensions_heterogeneous}  In light of the trace theorem (see~Theorem~\ref{thm:trace}), it is natural to expect that functions in $\Hr$ can be extended to $\R^n$ in such a way that their restriction to the complement of $\overline{\Omega}$ lies in $W^{1,p}(\overline{\Omega}^c)$. As we will confirm, these extensions lie in the following space. 

\begin{definition}[Heterogeneous nonlocal Sobolev spaces on $\R^n$]\label{def:newspacesRn}
With $\Grr=\Qr\nabla$ the extended heterogeneous nonlocal gradient from  Definition~\ref{def:hetnonlocgrad2}, we define
\[
\Hrr:=\{u \in L^p(\R^n) \,:\, \Grr u \in L^p(\R^n;\R^n)\},
\]
endowed with the norm  
\[
\norm{u}_{\Hrr}:=\norm{u}_{L^p(\R^n)} +\norm{\Grr u}_{L^p(\R^n;\R^n)}
\]
for $u\in \Hrr$.
\end{definition}

Let us now discuss some properties of the function spaces $\Hrr$. Since they are natural analogues of $\Hr$, the statements are essentially the same and the proofs are parallel if one replaces in the arguments the restricted auxiliary operators of Section~\ref{sec:operators}\,a) and d), that is $\PROmega, \RrOmega, \Cr^{\Omega}$, by the corresponding 
 operators $\PR, \Rr, \Cr$.  
We collect here the aforementioned statements 
 for easier reference. 
In fact, the spaces $\Hrr$ are reflexive, separable Banach spaces for $p\in (1,\infty)$ satisfying the continuous embeddings
\begin{align*}
W^{1,p}(\R^n)\hookrightarrow \Hrr \hookrightarrow H^{\lambda, p}(\R^n),
\end{align*}
cf.~Proposition~\ref{prop:clasembedding}.
The next result is the counterpart of Lemma~\ref{le:connection}.

\begin{lemma}[Translation mechanism in the extended setting]\label{lem:connectionRn}  These statements hold:
\begin{itemize}
\item[$(i)$] The operator $\Qr:\Hrr \to W^{1,p}(\R^n)$ is bounded and
\begin{align*}
\nabla \Qr - \Grr \in \LBH^+. 
\end{align*}
\item[$(ii)$] The operator $\PR:W^{1,p}(\R^n) \to \Hrr$ is bounded and
\begin{align*}
\Grr \PR - \nabla \in \LBH^\gamma \quad\text{for all $\gamma>1-\lambda$. } 
\end{align*}
\end{itemize}
\end{lemma}
In analogy to Corollary~\ref{cor:represHrhop}, we obtain that
\begin{align}\label{HrrW1p}
\Hrr=\{u\in L^p(\R^n): \Qr u\in W^{1,p}(\R^n)\},
\end{align}
with the norm  $\norm{\,\cdot\, }_{\Hrr}$ equivalent to $\norm{\Qr \,\cdot\,}_{W^{1,p}(\R^n)}$. 
In particular, it holds for any $u \in L^p(\R^n)$ that 
$\Qr u \in W^{1,p}(\R^n)$ if and only if $u \in \Hrr$. 
This observation can be used to deduce the following representation result for $\Hrr$. 

\begin{corollary}[Representation of $\Hrr$]\label{prop:hrrproperties}
It holds that 
\begin{align}\label{X}
\Hrr= \bigl\{u \in L^p(\R^n) \,:\, u|_{\Omega} \in \Hr, \ u|_{\overline{\Omega}^c}\in W^{1,p}(\overline{\Omega}^c), \ \Trho(u|_\Omega)= T_{W^{1,p}(\overline{\Omega}^c)}(u|_{\overline{\Omega}^c})\bigr\},
\end{align}
and $\norm{u}_{\Hrr} = \norm{u|_{\Omega}}_{\Hr} + \norm{u|_{\overline{\Omega}^c}}_{W^{1,p}(\overline{\Omega}^c)}$ for $u\in \Hrr$. 
\end{corollary}
\color{black}
\begin{proof}
We start by observing that
\begin{align}\label{Qr}
\Qr u =\begin{cases}
\Qro (u|_\Omega) &\text{on $\Omega$},\\
u|_{\overline{\Omega}^c} &\text{on $\overline{\Omega}^c$},
\end{cases} \qquad \text{for any $u\in L^p(\R^n)$,}
\end{align}
due to~\eqref{Qcal_restriction} and~\eqref{dualQast}, cf.~also~\eqref{eq:Qrrepr}. \color{black} Let $X$ denote the function space on the right-hand side of~\eqref{X}. 

If $u \in X$, then $\Qr u$ 
lies in $W^{1,p}(\R^n)$. Indeed, in view of~\eqref{Qr} it is the concatenation of two Sobolev functions with matching trace values on $\partial\Omega$, see also~\eqref{def:traceoperator}. Hence,~\eqref{HrrW1p} yields $u \in \Hrr$.

Conversely, if $u\in \Hrr$, it is clear that
\begin{equation}\label{eq:extensiongradient}
\Grr u=\begin{cases}
\Gr(u|_\Omega) & \text{on $\Omega$},\\
\nabla u|_{\overline{\Omega}^c} & \text{on $\overline{\Omega}^c$,}
\end{cases}
\end{equation}
which implies $u|_{\Omega} \in \Hr$ and  $u|_{\overline{\Omega}^c} \in W^{1,p}(\overline{\Omega}^c)$.
Moreover, since $\Qr u \in W^{1,p}(\R^n)$ by~\eqref{HrrW1p} with $\Qr u|_{\overline{\Omega^c}}=u|_{\overline{\Omega^c}}$, the compatibility of traces follows given the definition of the trace operator on $\Hr$ in~\eqref{def:traceoperator}, so that $u\in X$.

The representation formula for the norm $\norm{\cdot}_{\Hrr}$ is immediate to see with~\eqref{eq:extensiongradient}.
\end{proof}


With our improved understanding of $\Hrr$, we will now show that the extension operator of Section~\ref{sec:operators}\,c) also provides the desired extension of functions in the heterogeneous nonlocal Sobolev spaces. 
\begin{corollary}[Extension operator]\label{cor:extension}
The operator $\Er$ introduced in~\eqref{eq:extension} maps continuously from $\Hr$ to $\Hrr$ and satisfies $\Er u|_{\Omega}=u$ for all $u \in \Hr$. 
\end{corollary}
\begin{proof}
The fact that $\Er u|_{\Omega}=u$ for all $u \in \Hr$ is immediate from the definition. For the boundedness, we use~\eqref{QEswitch}, Corollary~\ref{cor:represHrhop} and \eqref{HrrW1p} to find that
\begin{align*}
\norm{\Er u}_{\Hrr} &\leq C\norm{\Qr \Er u}_{W^{1,p}(\R^n)} = C\norm{E \Qro u}_{W^{1,p}(\Omega)} \\ &\leq C\norm{\Qro u}_{W^{1,p}(\Omega)} \leq C\norm{u}_{\Hr}. \qedhere
\end{align*}
\end{proof}

 \subsection{Regularity results}
We have previously seen that $W^{1,p}(\Omega)$ is contained in $\Hr$, and that elements of $\Hr$ have the same boundary traces as functions in $W^{1,p}(\Omega)$. Next, we will show that, away from the boundary of $\Omega$, the functions in $\Hr$ can exhibit lower regularity than classical Sobolev functions. 
 
 To make this precise, let us recall that $H^{\kappa,p}_0(U)$ with $U\subset \R^n$ open and bounded is the closed subspace of functions in the Bessel potential space $H^{\kappa, p}(\Rn)$ that are zero in the complement of $U$, see~\eqref{Besselzero}, with $\kappa\in (0,1)$ the parameter in~\ref{itm:upper}. 

\begin{lemma}[Refined embedding into $\Hr$]
\label{le:discontinuous}
For any open set $U \Subset \Omega$, it holds that $H^{\kappa,p}_0(U) \hookrightarrow \Hrr$. 
In particular, also $H^{\kappa,p}_0(U)|_{\Omega} \hookrightarrow \Hr$. 
\end{lemma}

\begin{proof}
 The proof is based on the idea that neglecting boundary effects allows to relate $\Gr$ with a pseudo-differential operator associated to a cut-off version of the symbol $\qr$, which can be shown to have a lower order, namely $\kappa-1$ instead of $0$.

{We assume that $U$ is a Lipschitz domain, since we can always replace $U$ by a larger set compactly contained in $\Omega$.} Given the representation of $\Qr$ in \eqref{eq:Qrrepr}, one can find a set $U' \Subset \Omega$ such that $\supp(\Qr \varphi) \subset U'$ for all $\varphi \in C_c^\infty(U)$. By taking a cut-off function $\chi \in C_c^{\infty}(\Omega)$ with $\chi \equiv 1$ on $U'$, it follows that
\begin{equation}\label{eq:qrcutoff}
\Qr \varphi= \chi \Qr \varphi = \Op(\qr^{\chi}) \varphi \quad \text{for all $\varphi \in C^\infty_c(U)$},
\end{equation}
with the symbol 
\begin{align*}
\qr^{\chi}(x,\xi):=\chi(x)\qr(x,\xi)=\chi(x)\widehat{Q}_{\rho_1}(\d(x)\xi) \qquad\text{  for $x,\xi \in \R^n$. }
\end{align*} 
We will  verify that $\qr^{\chi}$ lies in the H\"ormander class of order $\kappa-1$. 
Indeed, since \begin{align*}
\abs{\widehat{Q}_{\rho_1}} \leq C \abss{\cdot}^{\kappa-1}
\end{align*}  
 by Proposition~\ref{prop:BeMCSc24} and \ref{itm:upper} and since the horizon function $\d$ is strictly positive on the support of $\chi$, say $\delta\geq \delta_\chi>0$ on $\supp \chi$, we may use the estimate in Lemma~\ref{le:symbolbounds} to conclude 
for any $\mu\in (0,1)$ and all $\alpha,\beta \in \N_0^{n}$ that
\begin{align*}
\abs{\partial^\beta_x \partial^{\alpha}_\xi \qr^\chi(x,\xi)} & \leq C_{\alpha,\beta} \abss{\xi}^{-\abs{\alpha}+\mu\abs{\beta}}|\chi(x)|\,\abs{\widehat{Q}_{\rho_1}(\d(x)\xi)} 
\\ &\leq C_{\alpha, \beta}   \abss{\xi}^{ -\abs{\alpha}+\mu\abs{\beta}} \abss{\delta_\chi \xi}^{\kappa-1}  \leq C_{\alpha, \beta}   \abss{\xi}^{\kappa -1 -\abs{\alpha}+\mu\abs{\beta}}
\end{align*}
for all $x,\xi\in \R^n$. 
This yields $\qr^{\chi} \in S^{\kappa-1}_{1,\mu}$ for any $\mu\in (0,1)$.

Exploiting the mapping properties of the associated pseudo-differential operator $\Op(\qr^{\chi})$ along with the density of $C_c^\infty(U)$ in $H^{\kappa, p}_0(U)$ implies by~\eqref{eq:qrcutoff} that
\[
\norm{\Qr u}_{W^{1,p}(\Rn)} =\norm{\Op(\qr^{\chi}) u}_{W^{1,p}(\Rn)} \leq C\norm{u}_{H^{\kappa,p}(\R^n)} \quad \text{for all $u \in H^{\kappa,p}_0(U)$}. 
\]
Finally, we exploit~\eqref{HrrW1p} together with the equivalence of norms stated below it to conclude the proof of the first embedding.
The second one follows simply by restriction of these spaces to $\Omega$. 
\end{proof}

We can now combine the previous lemma with Proposition~\ref{prop:clasembedding} to formulate our next result, which summarizes the regularity properties of $\Hr$ at a glance.
Roughly speaking, $\Hr$ lies in between the Bessel potential spaces with the fractional orders $\kappa$ from~\ref{itm:upper} and $\lambda$~from~\ref{itm:lower},  away from the boundary.

\begin{corollary}[Regularity properties of $\Hr$]\label{cor:hrregularity}
For any open set $U \Subset \Omega$, it holds that
\[
H^{\kappa,p}_0(U)|_{\Omega} + W^{1,p}(\O) \subset \Hr \subset H^{\lambda,p}(\O).
\]
\end{corollary}
Considering the properties of Bessel potential spaces, it is immediate to conclude that $\Hr$ may contain also less regular functions, for instance with discontinuities in the form of cavitations or jumps.

\begin{example}\label{ex:jumps}
a) If $\kappa p < n$, the fact that there exist discontinuous functions in $H^{\kappa,p}_0(B_1(0))$ (see e.g.~\cite[Lemma~2.5]{BeCuMC} for the example $x\mapsto x/|x|\chi(|x|)$ with a suitable smooth cut-off function $\chi$), allows us to infer the appearance of discontinuities in elements of $\Hr$ as well. 
In the case $\kappa p <1$, it even holds that {$\mathbbm{1}_{(0,1/n)^n} \in H^{\kappa,p}_0(B_1(0))$} (cf.~\cite[Lemma~2.4]{BeCuMC}), illustrating that $\Hr$ contains functions admitting discontinuities across hypersurfaces. \smallskip

\color{black} b) Consider an admissible kernel function $\rho_1$ satisfying~\ref{itm:lower} and~\ref{itm:upper} with $\lambda=\kappa=s\in (0,1)$, such as the truncated fractional kernel in Example~\ref{ex:rho1}. Then, Corollary~\ref{cor:hrregularity} implies 
\begin{align*}
\Hr|_U=H^{s, p}(U) \quad \text{ for any open $U\Subset \Omega$, }
\end{align*}
which characterizes the space $\Hr$ almost completely,  up to the behavior near the boundary $\partial \Omega$.
\end{example}

\subsection{Poincar\'{e} inequalities and functions with zero nonlocal gradient} Let us now turn to proving a Poincar\'{e} inequality for the gradient $\Gr$. We will first establish the following abstract version, which follows from a standard contradiction argument.
Functions with vanishing heterogeneous nonlocal gradient will play a central role in the analysis of this section, and are also of independent interest.  We denote the closed subspace of $\Hr$ corresponding to the kernel of $\Gr$ as 
\[
\Nr:=\bigl\{u \in \Hr \,:\, \Gr u = 0 \ \text{a.e.~in $\Omega$}\bigr\}. 
\]

\begin{lemma}\label{le:abstractpoincare}
Let $X \subset \Hr$ be a closed subset such that $X \cap \Nr = \{0\}$. Then, there is a constant $C>0$ such that
\[
\norm{u}_{L^p(\Omega)} \leq C \norm{\Gr u}_{L^p(\Omega;\R^n)} \quad\text{for all $u \in X$}.
\]
\end{lemma}

\begin{proof}
Suppose to the contrary that there exists a sequence $(u_j)_j \subset X$ with
\[
1=\norm{u_j}_{L^p(\Omega)}  > j \norm{\Gr u_j}_{L^p(\Omega;\R^n)}
\]
for all $j \in \N$.  This implies in particular that, up to selecting a non-relabeled subsequence, $u_j \weakto u$ in $\Hr$ as $j\to \infty$ for some $u \in \Hr$. In fact, since $\Gr u_j \to 0$ in $L^p(\Omega;\R^n)$, we deduce that $\Gr u = 0$, and thus, $u\in \Nr$. Proposition~\ref{prop:clasembedding} along with the compact embedding $H^{\lambda,p}(\Omega) \hookrightarrow\hookrightarrow L^p(\Omega)$, which is due to~\eqref{Hcompactembedding}, implies even the strong convergence $u_j \to u$ in $\Hr$. 
\color{black}The closedness of $X$ then yields that $u \in X$, from which we conclude that $u=0$, given that $X\cap \Nr=\{0\}$ by assumption. Consequently, 
\[
1=\lim_{j \to \infty}\norm{u_j}_{L^p(\Omega)}=\norm{u}_{L^p(\Omega)} = 0
\]
produces a contradiction, as desired. 
\end{proof}

We point out that the result  of Lemma~\ref{le:abstractpoincare} remains somewhat implicit as long as $\Nr$ is not precisely characterized. This motivates us to work toward a better understanding of functions with vanishing nonlocal gradients $\Gr$ and ideally, obtain an explicit description of $\Nr$. While in the classical case the analogue problem is easily solved, as the kernel consists solely of all constant functions, it is more subtle in nonlocal settings. For constant-horizon nonlocal gradients, the authors recently showed in~\cite{KrS24} that the corresponding kernel forms an infinite-dimensional vector space, with each function identified uniquely by its boundary values in a collar region and a mean value condition. 

In  contrast, in our set-up of nonlocal gradients with a space-dependent horizon that decays sufficiently fast toward the boundary, it turns out that $\Nr$ is finite-dimensional. Under an additional assumption -- specifically, that the horizon function satisfies a property we refer to as mild variation, defined below -- we can prove more and identify the space with the constant functions, as in the classical case. 

Before we state these findings in Theorem~\ref{th:Nchar} below, let us first introduce what we mean by a mildly varying horizon function: 
Consider $ \delta\in \Cinfo$ as in Section~\ref{subsec:defnonlocalgradient}, satisfying $(i)$-$(iii)$ of Lemma~\ref{le:constructioneta}, and rewrite it in the form 
\begin{align*}
\delta(x) =\bar{\delta}\eta_\delta(x) \quad \text{for $x\in \R^n$,}
\end{align*} where $\bar{\delta}>0$ denotes the maximal value of $\delta$ on $\R^n$ and $ \eta_\delta \in  \Cinfo$ is the normalization of $\delta$ with maximum $1$; note that $\eta_\delta$ is again an admissible horizon function, with the associated kernel function $\rho_{\eta_\delta(\cdot)}$. The symbols defined according to \eqref{symbol_qrho} have the relation
\begin{align*}
q_{\rho(\cdot)}(x,\xi)=  \widehat{Q}_{\rho_1}(\bar{\delta} \eta_\delta(x)\xi) = q_{\rho_{\eta_\delta}(\cdot)}(x, \bar{\delta}\xi)\quad \text{for $x, \xi\in \R^n$.}
\end{align*}

In view of Lemma~\ref{le:symbolbounds}, we find that $q_{\rho_{\eta_\delta}(\cdot)}\in S^{0}_{1, \mu}$ for $\mu\in (0,1)$ fulfills the requirements of Proposition~\ref{prop:exactinvertible} with $m'=\lambda-1$ as long as we take $\mu < \lambda$. 
{Hence, by Proposition~\ref{prop:exactinvertible}, there exists for each $p \in (1,\infty)$ a $\bar{\d}_0 >0$ such that the kernel of $\Qr:L^p(\R^n)\to L^p(\R^n)$ is trivial if $ \bar{\delta} < \bar{\d}_0$.} The condition to be introduced is intended to  guarantee  injectivity of the operator $\Qr$: 
\begin{align}\label{def:mildlyvarying}
\text{ We call the horizon function $\delta$ \textit{mildly varying} if $\bar{\delta}<\bar \delta_0$ holds. }
\end{align}

\color{black}
After these preparations, we are now ready to formulate our result on the kernel $\Nr$.

\begin{theorem}[Properties of the kernel $\Nr$]\label{th:Nchar}
It holds that $ \Nr\subset \Cinf$ is finite-dimensional. If the horizon function $\delta$ is mildly varying, then 
\begin{align}\label{Nrconstant}
\Nr = \{ u: \Omega\to \R\,:\, u \ \text{is constant}\}.
\end{align}
\end{theorem}

\begin{proof}
Let $u \in \Nr$, and define $\tilde{u} = E_{\rho(\cdot)}u \in \Hrr \subset H^{\lambda,p}(\R^n)$. The proof relies on a comparison of the heterogeneous nonlocal gradient of the extension $\tilde u$ with its classical (distributional) gradient. 
To this end, we introduce the auxiliary function
\[
\tilde w:=\nabla\tilde{u}-\Grr \tilde{u} \in H^{\lambda-1,p}(\R^n;\R^n),
\]
observing that $\tilde w|_\Omega=\nabla u$  as well as $\tilde w|_{\overline{\Omega}^c}=0$,  in view of \eqref{eq:extensiongradient} along with $\Gr u =0$. 

First, we will
establish that 
\begin{align}\label{Qtildew=0}
\Qr \tilde w=0.
\end{align} 
As $\Grr \tilde{u} \in L^p(\R^n;\R^n)$ and $\Qr \in \LBH^0$, it follows that
\[
\Qr \tilde w = \Grr \tilde{u}-\Qr \Grr \tilde{u} \in L^p(\R^n;\R^n).
\]
For the restrictions to $\Omega$ and $\overline{\Omega}^c$, we infer from the relation between $\Grr$ and $\Gr$ in \eqref{eq:extensiongradient} that
\[
(\Qr \tilde w)|_{\Omega}=\Gr u - \Qro \Gr u = 0\qquad  \text{and }\qquad (\Qr \tilde w)|_{\overline{\Omega}^c} =\tilde w|_{\overline{\Omega}^c}= \nabla \tilde{u} - \nabla \tilde{u} = 0,
\]
respectively.
Considering $\Qr \tilde w \in L^p(\R^n;\R^n)$, this yields~\eqref{Qtildew=0}. 

We then find in light of \eqref{eq:res} that
\begin{align}\label{wineigenspace}
\tilde w = \PR \Qr \tilde w - \Rr \tilde w = -\Rr \tilde w \  \in C^{\infty}(\R^n;\R^n), 
\end{align}
which implies, together with \eqref{Rprop_complementary}, that $\tilde w \in C_0^{\infty}(\overline{\Omega};\R^n)$. 
Since $\nabla u=\tilde w|_\Omega$, we deduce $u\in \Cinf$, 
yielding the stated regularity of functions in $\Nr$. 

That $\Nr$ has finite dimension can be derived also from~\eqref{wineigenspace}. 
Indeed, the identity in~\eqref{wineigenspace} indicates that  $\tilde w$ lies in the eigenspace of $\Rr$ corresponding to the eigenvalue $-1$. By exploiting the compactness of $\Rr:L^p(\R^n)\to L^p(\R^n)$ as a smoothing operator in $\LBH^{-\infty}$ satisfying \eqref{Rprop_complementary}, we conclude from standard spectral theory, e.g.~\cite[Theorem~7.7.1]{Con07}, that this eigenspace is a finite-dimensional subspace of $C^\infty(\R^n;\R^n)$. Given $\nabla u=\tilde w|_\Omega$ for $u\in \Nr$, this implies that the kernel $\Nr$ can be at most finite-dimensional, as desired.

\color{black}
To prove the second part of the statement, we make use of the additional assumption that $\delta$ is mildly varying, which allows us to exploit that the kernel of $\Qr$ is trivial. From~\eqref{Qtildew=0}, we thus deduce $\tilde w=0$, and consequently, $\nabla u=w= \tilde w|_{\Omega}=0$. { Due to the connectedness of $\Omega$, it follows that $u$ is constant, which shows~\eqref{Nrconstant} and concludes the proof.}
\end{proof}

\begin{remark}
While our current proof does not extend to this case, we expect that \eqref{Nrconstant} holds also when $\delta$ is not mildly varying. In fact, this is corroborated by some numerical simulations, which show that the eigenvalues of the discretization of $\Qro$ are all positive and bounded away from zero.
\end{remark}

\color{black}
By combining Theorem~\ref{th:Nchar} with Lemma~\ref{le:abstractpoincare}, we obtain the following collection of Poincar\'e inequalities, which play a key role in the analysis of variational problems and partial differential equations involving heterogeneous nonlocal gradients. 
\begin{corollary}\label{cor:poincare}
Suppose the horizon function $\delta$ varies mildly. Then, there exists a constant $C>0$ such that
\[
\norm{u}_{L^p(\Omega)} \leq C \norm{\Gr u}_{L^p(\Omega;\R^n)} \quad\text{for all $u \in X$},
\]
with $X \subset \Hr$ as in any of these three settings: \smallskip
\begin{itemize}
\item[(i)] $X=\mathring{H}^{\rho(\cdot),p}(\Omega):=\left\{u \in \Hr \,:\, \int_\Omega u\,dx =0\right\}$; \smallskip
\item[(ii)] $X= \left\{u \in \Hr \,:\, \Trho u=0 \ \text{ $\Hcal^{n-1}$-a.e.~on $\Gamma$}\right\}$ for any { Borel-measurable} set $\Gamma \subset \partial \Omega$ with $\Hcal^{n-1}(\Gamma) >0$; \smallskip
\item[(iii)] $X= \left\{u \in \Hr \,:\, u=0 \ \text{a.e.~on $U$}\right\}$ for any { measurable} set $U \subset \Omega$ with $\abs{U}>0$.
\end{itemize}
\end{corollary}

\color{black}
\section{Existence theory for heterogeneous nonlocal variational problems}\label{sec:variational}

This section discusses applications of the previously established properties of heterogeneous nonlocal gradients in the context of nonlocal variational problems and their associated boundary value problems.
Our focus is on questions related to the existence of minimizers for integral functionals of the form
\begin{align*}
\Fr(u)=\int_{\Omega} f(x,u,\Gr u)\,dx  \quad \text{for $u \in \Hrm$},
\end{align*}
where $p\in (1, \infty)$, $\Omega\subset \R^n$ is a bounded { and connected} Lipschitz domain and $f:\Omega\times \R^N\times \R^{N \times n}\to \R_\infty:=\R\cup\{\infty\}$  a Carath\'{e}odory integrand with suitable growth and coercivity properties. Throughout, we assume that $\Gr$ is as introduced in Section~\ref{subsec:defnonlocalgradient}, with a horizon function $ \delta\in \Cinfo$ that satisfies $(i)$-$(iii)$ of Lemma~\ref{le:constructioneta}. Wherever  a Poincar\'e inequality is used, we will require in addition that $\delta$ is mildly varying, cf.~\eqref{def:mildlyvarying} and Corollary~\ref{cor:poincare}.

\subsection{Weak lower semicontinuity} A central ingredient in applying the direct method to the nonlocal variational problems under consideration is the weak lower semicontinuity of the functionals $\Fr$.

We begin with a characterization result in the case of integrands with standard $p$-growth. The following theorem shows that the quasiconvexity of $f$ (in the sense of Morrey~\cite{Morrey}) with respect to its third argument is both necessary and sufficient for the weak lower semicontinuity of $\Fr$ -- exactly as in the well-known case of local integral functionals depending on classical gradients, e.g.~\cite{AcF84,Morrey}.  For analogous results for the Riesz fractional gradient and nonlocal gradients with homogeneous horizon, see~\cite{KrS24, CuKrSc23}.  
The proof strategy builds on the connection between the classical and nonlocal frameworks, established through the (non-perfect) translation mechanism of Lemma~\ref{le:connection}, which enables us to reduce the analysis to the classical case. Specifically, the mismatch between $\Gr$ and the translated classical gradient involves terms with lower-order operators, which can be controlled using strong convergence arguments. 

\begin{theorem}[Characterizing weak lower semicontinuity via quasiconvexity]\label{th:lscchar}
Let $f:\Omega\times \R^N\times \Rmn \to \R$ be a Carath\'{e}odory integrand that satisfies 
\begin{align}\label{eq:growth}
 -C(1+\abs{z}^p+\abs{A}^q)  \leq f(x,z,A) \leq C(1+\abs{z}^p+\abs{A}^p) 
\end{align}
for a.e.~$x \in \Omega$ and all $(z,A) \in \R^N\times\Rmn$, with a constant $C>0$ and some $q \in (1,p)$.
Then, the functional $\Fr:\Hrm \to \R$ given by
\[
\Fr(u)=\int_{\Omega} f(x,u,\Gr u)\,dx 
\]
is weakly lower semicontinuous on $\Hrm$ if and only if $f(x,z,\cdot)$ is quasiconvex, 
 i.e., 
\begin{align*}
f(x,z,A)\leq \inf_{v\in W^{1, \infty}_0((0,1)^n;\R^N)}\int_{(0,1)^n} f(x, z, A+\nabla v(y)) \, dy \qquad \text{for every $A\in \R^{N \times n}$,}
\end{align*} 
 for a.e.~$x \in \Omega$ and all $z \in \R^N$.
\end{theorem}
\begin{proof}
\textit{Sufficiency.} Let $(u_j)_j \subset \Hrm$ converge weakly to some $u \in \Hrm$ as $j\to \infty$, that is, $u_j\weakly u$ in $\Hrm$. Lemma~\ref{le:connection} then implies that
\begin{align*}
v_j:=\Qro u_j \weakto \Qro u =: v \quad \text{ in $W^{1,p}(\Omega;\R^N)$.}
\end{align*}
 Moreover, by Proposition~\ref{prop:clasembedding} and the compact embedding $H^{\lambda,p}(\Omega;\R^N)\hookrightarrow\hookrightarrow H^{\mu,p}(\Omega;\R^N)$ for any $\mu \in (0, \lambda)$ (cf.~Section~\ref{sec:bps}), we deduce that $u_j \to u$ in $H^{\mu,p}(\Omega;\R^N)$. Using the commutator from \eqref{eq:cro} with the mapping properties~\eqref{mappingprop_C}, it follows that $\Cro u_j \to \Cro u$ in $L^p(\Omega;\Rmn)$, and hence, 
\begin{equation}\label{eq:residual}
R_j:=\Gr u_j - \nabla v_j = -\Cro u_j \to -\Cro u = \Gr u- \nabla v=:R \quad \text{in $L^p(\Omega;\Rmn)$}.
\end{equation}

The rest of the argument relies on common techniques from Young measure theory; for introductions into this topic, we refer~e.g.~to~\cite{Rindler, FoL07}. Let $\nu=(\nu_x)_{x\in \Omega}$ be the gradient Young measure generated, up to a subsequence, by the sequence $(\nabla v_j)_j$.  In view of \eqref{eq:residual} in combination with~\cite[Corollary 8.7]{FoL07}, one obtains that $(\Gr u_j)_j$ generates, up to a subsequence, the Young measure $\nu$ translated by $R$.
Since $u_j \to u$ in $L^p(\O;\R^N)$, we know that $(u_j)_j$,  again up to selecting a subsequence, generates the Dirac measures $(\delta_{u(x)})_{x\in \Omega}$ as Young measures (cf.~\cite[Lemma~4.12]{Rindler}).

Based on these observations, the fundamental theorem for Young measures (see~e.g.~\cite[Theorem~8.6\,$(i)$]{FoL07}) is applicable here under consideration of the lower bound in \eqref{eq:growth}, along with the result on pairs of sequences in~\cite[Lemma~5.19]{Rindler}. This allows us to conclude
\begin{align*}
\liminf_{j \to \infty} \Fr(u_j) &=\liminf_{j \to \infty} \int_\Omega f(x,u_j(x),\Gr u_j(x))\,dx \\
&\geq \int_{\Omega} \int_{\Rmn} f(x,u(x),A+R(x)) \,d\nu_x(A)\,dx \\
&\geq \int_{\Omega} f(x,u(x), \nabla v(x) + R(x))\,dx\\
&=\int_{\Omega} f(x,u(x),\Gr u(x))\,dx =\Fr(u),
\end{align*} 
where the second estimate corresponds to the generalized Jensen's inequality applied to the quasi\-convex function $A \mapsto f(x,u(x), A+R(x))$ for a.e.~$x \in \Omega$ (cf.~\cite[Theorem~7.15]{Rindler}). \smallskip

\textit{Necessity.} Let $(x_0,z_0,A_0) \in \R^n\times\R^N \times \Rmn$ and $(v_j)_j \subset W^{1,\infty}(\Omega;\R^N)$ be such that $v_j \starto 0$ in $W^{1,\infty}(\Omega;\R^N)$ as $j\to \infty$. Our aim is to show that
\begin{align}\label{eq:goal}
\int_{\Omega}f(x,z_0+A_0(x-x_0), A_0) \,dx \leq \liminf_{j \to \infty} \int_{\Omega} f(x,z_0 + A_0(x-x_0), A_0+\nabla v_j)\,dx,
\end{align}
from which the quasiconvexity of $f$ follows according to~\cite[Lemma~3.18 and Remark~3.19\,(ii)]{Dac08}. 

In the following, set  $u_0(x):=z_0+A_0(x-x_0)$ for $x \in \Omega$ and note that
\begin{align}\label{affine_gradient}
 \Gr u_0(x)= (\Qr^\Omega \nabla u_0) (x) = \int_{\R^n} \d(x)^{-n}Q_{\rho_1}\left(\frac{x-y}{\d(x)}\right)A_0\,dy= \norm{Q_{\rho_1}}_{L^1(\R^n)}A_0 = A_0,
\end{align} 
for all $x \in \Omega$, cf.~\eqref{eq:Qrrepr} and~\eqref{eq:rhonormalization}. Further, we define 
 \begin{align*}
 u_j:=\PROmega v_j \in H^{\rho(\cdot),p}(\Omega;\R^N) \qquad \text{for $j \in \N$}
 \end{align*} 
and observe that $u_j\weakly 0$ in $H^{\rho(\cdot), p}(\Omega;\R^N)$ as $j\to \infty$ by Proposition~\ref{le:connection}\,$(ii)$. Combining the latter  with Proposition~\ref{prop:clasembedding} and~\eqref{Hcompactembedding}, shows that $\PROmega$ maps compactly from $W^{1, p}(\Omega;\R^N)$ into $L^p(\Omega;\R^N)$, and we obtain that $\Gr\PROmega - \nabla$ is a compact  operator $W^{1, p}(\Omega;\R^N)\to L^p(\Omega;\R^N)$ in view of~\eqref{translation2} applied with some $\gamma\in (1-\lambda, 1)$ and~\eqref{Hcompactembedding}. 
Hence, 
 \begin{align}\label{ujtozerostrong}
 u_j \to 0 \quad \text{in $L^p(\Omega;\R^N)$}
   \end{align}
   and
\begin{equation}\label{eq:residualzero}
R_j:= \Gr u_j-\nabla v_j \to 0 \quad \text{in $L^p(\Omega;\Rmn)$}.
\end{equation}

We can then exploit the weak lower semicontinuity of $\Fr$ to obtain
\begin{align*}
\int_{\Omega}f(x,u_0, \Gr u_0) \,dx &\leq \liminf_{j \to \infty} \int_\Omega f(x,u_0+u_j, \Gr u_0+\Gr u_j)\,dx\\
& = \liminf_{j \to \infty}\int_\Omega f(x,u_0+u_j, \Gr u_0+\nabla v_j+R_j)\,dx\\
&=\liminf_{j \to \infty} \int_\Omega f(x,u_0, \Gr u_0+\nabla v_j)\,dx;
\end{align*}
note that \eqref{ujtozerostrong} and \eqref{eq:residualzero} are used in the last line, as well as
 the $L^p$-equi-integrability of $(\nabla v_j)_j$ as a bounded sequence in $W^{1, \infty}(\Omega;\R^N)$ (see~e.g.~\cite[Lemma~4.10]{KrS22}). In light of the definition of $u_0$ and~\eqref{affine_gradient}, this shows \eqref{eq:goal}, as desired.
\end{proof}

For extended-valued integrands that do not satisfy the upper bound in \eqref{eq:growth}, the sufficiency part of Theorem~\ref{th:lscchar} still holds, provided $f$ falls witin the more restrictive setting of polyconvexity as in \cite[Theorem~8.16]{Dac08}; recall that a function is polyconvex if it can be expressed as a convex functions of the minors of its argument, cf.~\cite{Ball1977}. This is particularly relevant in applications of hyperelasticity, where such integrand functions are used as energy densities, see e.g.~\cite[Example~6.2-6.4]{Rindler}. 

\begin{theorem}[Polyconvexity as sufficient condition]\label{th:lscsuf}
Let $p>\min\{N,n\}$ and $f:\Omega\times \R^N\times \Rmn \to \R_\infty$ be a Carath\'{e}odory integrand satisfying the lower bound in~\eqref{eq:growth}. If $f(x,z, \cdot)$ is polyconvex\footnote{More precisely, $f(x,z,A)=\tilde{f}(x,z,T(A))$ for all $(x,z,A) \in \Omega \times \R^N \times \Rmn$, where $T(A)$ is the vector of minors of $A$ and $\tilde{f}$ is a Carath\'eodory function that is convex in the third argument and satisfies an analogous lower bound to \eqref{eq:growth}.}  for a.e.~$x \in \Omega$ and all $z \in \R^N$, then the functional $\Fr:\Hrm \to \R_\infty$ given by
\[
\Fr(u):=\int_{\Omega} f(x,u,\Gr u)\,dx 
\] 
is weakly lower semicontinuous on $\Hrm$. 
\end{theorem}

\begin{proof}
It suffices to show that the minors of the heterogeneous nonlocal gradient are weakly continuous, as we do below. With this observation at hand, the rest of the argument is standard in the literature, for instance, one can follow~\cite[Theorem~8.16]{Dac08}.

Let $M_k:\Rmn \to \R$ be any minor of order $k \leq \min\{N,n\}$. For any $a \in L^{\infty}(\O)$, the auxiliary function  
$$g:\Omega \times \Rmn \to \R, \quad g(x,A)=a(x)M_k(A),$$
as well as $-g$, satisfy both the assumptions of Theorem~\ref{th:lscchar}, given the fact that all minors are quasiaffine (see~e.g.~\cite[Theorem~5.20]{Dac08}) and the assumption $p>\min\{N,n\}$, which guarantees \eqref{eq:growth}. Hence, we deduce for any sequence $(u_j)_j$ with $u_j \weakto u$ in $\Hrm$ that
\[
\lim_{j \to \infty} \int_{\Omega} a(x) M_k(\Gr u_j)\,dx =\int_{\Omega} a(x) M_k(\Gr u)\,dx.
\]
Since this implies the weak convergence $M_k(\Gr u_j) \weakto M_k(\Gr u)$ in $L^{p/k}(\Omega)$, the proof is complete.
\end{proof}

\subsection{Existence results for minimizers}
By combining the previous lower semicontinuity statements with the Poincar\'e inequalities from Corollary~\ref{cor:poincare}, we obtain the following existence results under Dirichlet or natural boundary conditions via the direct method. 
Here, we restrict to the case with explicit, linear dependence of the integrand function on $u$, which is very common in models of elasticity for representing external forces. This assumption also leads to a simpler and more natural presentation of the associated Euler-Lagrange equations in the next subsection.   
We note, however, that the same proofs remain valid in the more general settings of Theorems~\ref{th:lscchar} and~\ref{th:lscsuf}.

\begin{corollary}[Existence of minimizers]\label{cor:existence}
Let $f:\Omega\times\Rmn \to \R_\infty$ be a Carath\'{e}odory integrand that satisfies  
\begin{align}\label{lowerboundf}
f(x, A) \geq c\abs{A}^p - C  \quad \text{for a.e.~$x \in \Omega$ and all 
$A \in \Rmn$},
\end{align}
with constants $C,c>0$. Suppose further that 
\begin{align*}
\text{  $f(x, \cdot)$ is quasiconvex and $f(x, \cdot)\leq C(1+|\cdot|^p)$}
  \end{align*}  
or that 
\begin{align*}
\text{$f(x, \cdot)$ is polyconvex 
and $p>\min\{N, n\}$}
\end{align*}
for a.e.~$x \in \Omega$. Then, if the horizon function $\delta$ is mildly varying, see \eqref{def:mildlyvarying}, the following holds: \smallskip
\begin{itemize}
\item[(i)] For any given $F \in L^{p'}(\Omega;\R^N)$ and $g \in W^{1-1/p,p}(\partial \Omega;\R^N)$, the functional
\[
u \mapsto \int_{\Omega} f(x,\Gr u)\,dx-\int_{\Omega}F\cdot u\,dx,
\]
admits a minimizer over $\Hrmg$, cf.~Definition~\ref{def:trace}. \smallskip

\item[(ii)] For any given $F \in L^{p'}(\Omega;\R^N)$ and $h \in L^{p'}(\partial \Omega;\R^N)$, the functional
\[
u \mapsto \int_{\Omega} f(x,\Gr u)\,dx-\int_{\Omega}F\cdot u\,dx-\int_{\partial\Omega} h\cdot \Trho u\,d\Hcal^{n-1},
\]
admits a minimizer over $\mathring{H}^{\rho(\cdot),p}(\Omega;\R^N)$, where $\Trho$ is the trace operator of Theorem~\ref{thm:trace}. 
\end{itemize}
\end{corollary}
\begin{proof}
This is an immediate application of the direct method in the calculus of variations. Indeed, the weak lower semicontinuity of the functionals is due to Theorems~\ref{th:lscchar} and~\ref{th:lscsuf}, in combination with the fact that the bounded linear functionals $u \mapsto \int_{\Omega}F\cdot u\,dx$ and $u \mapsto \int_{\partial\Omega} h\cdot T_{\rho(\cdot)} u\,d\Hcal^{n-1}$ are weakly continuous (cf.~Theorem~\ref{thm:trace}). The coercivity of the functionals follows by standard arguments from the lower bound on $f$ in~\eqref{lowerboundf} and the Poincar\'{e} inequalities in Corollary~\ref{cor:poincare}\,$(i)$ and $(ii)$, respectively.
\end{proof}

In a similar manner, one can obtain an existence result in a setting that combines the different boundary condition assumptions of $(i)$ and $(ii)$ in the previous corollary.
\begin{remark}[Mixed boundary conditions]\label{rem:mixed}
For a { Borel-measurable} $\Gamma \subset \partial \Omega$ with $\Hcal^{n-1}(\Gamma)>0$, one can deduce that the functional
\[
u \mapsto \int_{\Omega} f(x,\Gr u)\,dx-\int_{\Omega}F\cdot u\,dx-\int_{\partial\Omega \setminus \Gamma} h\cdot  T_{\rho(\cdot)} u\,d\Hcal^{n-1}
\]
with given functions $f, F, h$ as in~Corollary~\ref{cor:existence} has a minimizer over the subspace \begin{align*}
\{u \in \Hrm\,:\, T_{\rho(\cdot)} u=g \  \text{ $\Hcal^{n-1}$-a.e.~in $\Gamma$}\}.
\end{align*} This corresponds to mixed natural and Dirichlet boundary conditions, see~\eqref{eq:mixed} below.
\end{remark}

\subsection{Applications to nonlocal boundary value problems} 
Now that we have established the existence of minimizers for the functionals in Corollary~\ref{cor:existence}, we close by deriving the Euler-Lagrange equations that are satisfied by these minimizers. 

Suppose that $f$ is continuously differentiable in its second argument with
\[
\abs{D_Af(x,A)} \leq C(1+\abs{A}^{p-1}) \quad \text{for a.e.~$x \in \Omega$ and all $A \in \Rmn$},
\]
where $D_A f$ denotes the derivative of $f$ with respect to its second argument. 

It follows from standard methods that the minimizer $u \in \Hrmg$ of the functional in Corollary~\ref{cor:existence}\,$(i)$ satisfies
\[
\int_{\Omega} D_A f(x,\Gr u)\cdot \Gr \phi \,dx - \int_{\Omega}F\cdot \phi\,dx =0 \quad \text{for all $\phi \in C_c^{\infty}(\Omega;\R^N)$,}
\]
where we have used that $u+\phi \in \Hrmg$ for all $\phi \in C_c^{\infty}(\Omega;\R^N)$. In view of Proposition~\ref{prop:intbyparts}, the minimizer $u$ weakly solves 
\[
\begin{cases}
-\Div_{\rho(\cdot)}^{*}\left[D_A f(x,\Gr u)\right] = F & \text{in $\Omega$},\\
u=g &\text{on $\partial \Omega$}.
\end{cases}
\]
This is a  nonlocal system of partial differential equations subject to local Dirichlet boundary conditions.

Proceeding with the setting from Corollary~\ref{cor:existence}\,$(ii)$, we first assume the compatibility condition
\[
\int_{\Omega}F\,dx + \int_{\partial\Omega} h\,d\Hcal^{n-1} =0,
\]
in parallel to the classical case. Then, the minimizer $u \in \mathring{H}^{\rho(\cdot),p}(\Omega;\R^N)$ from Corollary~\ref{cor:existence}\,(ii) also minimizes the functional over the full space $\Hrm$. Hence, the usual techniques imply that $u$ satisfies
\[
\int_{\Omega} D_A f(x,\Gr u)\cdot \Gr \phi \,dx - \int_{\Omega}F\cdot\phi\,dx -\int_{\partial\Omega}h\cdot \phi\,d\Hcal^{n-1}=0 \quad \text{for all $\phi \in C^{\infty}(\R^n;\R^N)$.}
\]
Formally applying the integration by parts from Proposition~\ref{prop:intbyparts} and the fact that $\Gr = \nabla$ for $x \in \partial \Omega$, yields that $u$ is a weak solution of
\[
\begin{cases}
-\Div_{\rho(\cdot)}^{*}\left[D_A f(x,\Gr u)\right] = F & \text{in $\Omega$},\\
D_A f(x,\nabla u) \cdot \nu = h &\text{on $\partial\Omega$,}
\end{cases}
\]
with $\nu$ an outer unit normal to $\partial \Omega$. This corresponds to a system of nonlocal partial differential equations with classical natural boundary conditions on $\partial \Omega$.

Finally, we consider the setting from Remark~\ref{rem:mixed}, where the minimizer $u$ satisfies
\[
\int_{\Omega} D_A f(x,\Gr u)\cdot \Gr \phi \,dx - \int_{\Omega}F\cdot\phi\,dx -\int_{\partial\Omega\setminus \Gamma}h\cdot\phi\,d\Hcal^{n-1}=0,
\]
for all $\phi \in C_c^{\infty}(\R^n;\R^N)$ with $\phi=0$ on $\Gamma$. In the same way as before, we use the integration by parts from Proposition~\ref{prop:intbyparts} to deduce that $u$ weakly solves
\begin{equation}\label{eq:mixed}
\begin{cases}
-\Div_{\rho(\cdot)}^{*}\left[D_A f(x,\Gr u)\right] = F & \text{in $\Omega$},\\
D_A f(x,\nabla u) \cdot \nu = h &\text{on $\partial\Omega \setminus \Gamma$,}\\
u=g &\text{on $\Gamma$.}
\end{cases}
\end{equation}
Hence, we have the existence of weak solutions to this system of nonlocal partial differential equations with mixed local Dirichlet and natural boundary conditions.

\begin{example}
In the case $N=1$ and $f(x,A)=\frac{1}{2}\abs{A}^2$ for all $x\in \Omega$ and $A\in \R^{N \times n}$, the boundary value problem in \eqref{eq:mixed} reduces with $\Delta_{\rho(\cdot)}:=-\Div_{\rho(\cdot)}^{*}\Gr$ to 
\[
\begin{cases}
 \Delta_{\rho(\cdot)} u= F & \text{in $\Omega$},\\
\frac{\partial u}{\partial \nu} = h &\text{on $\partial\Omega \setminus \Gamma$,}\\
 u=g &\text{on $\Gamma$.}
\end{cases}
\]
This can be viewed as a nonlocal Laplace equation with mixed local boundary conditions.
\end{example}

  \section*{Acknowledgements}  
This research was funded in part by the Austrian Science Fund (FWF) projects \href{https://doi.org/10.55776/F65}{10.55776/F65} and \href{https://doi.org/10.55776/Y1292}{10.55776/Y1292}.   
Most of the research was done while H.S. was affiliated with
KU Eichst\"att-Ingolstadt and TU Wien, and H.S. also would like to thank the FNRS for its support in funding part of this work through a MIS-Ulysse project (Scientific Impulse Mandate Instrument)
number F.6002.25 during the finalizing stages of the paper. 
A visit of H.S. at KU Eichst\"att-Ingolstadt was supported by the Deutsche Forschungsgemeinschaft (DFG) through project 541520348, which the authors thankfully acknowledge.

\addcontentsline{toc}{section}{\protect\numberline{}References}
\bibliographystyle{abbrv}
\bibliography{bibliography}

\begin{thebibliography}{10}

\bibitem{AbP18}
H.~Abels and C.~Pfeuffer.
\newblock Characterization of non-smooth pseudodifferential operators.
\newblock {\em J. Fourier Anal. Appl.}, 24(2):371--415, 2018.

\bibitem{AcF84}
E.~Acerbi and N.~Fusco.
\newblock Semicontinuity problems in the calculus of variations.
\newblock {\em Arch. Rational Mech. Anal.}, 86(2):125--145, 1984.

\bibitem{ACFS25}
S.~Almi, M.~Caponi, M.~Friedrich, and F.~Solombrino.
\newblock A fractional approach to strain-gradient plasticity: beyond
  core-radius of discrete dislocations.
\newblock {\em Math. Ann.}, 391(3):4063--4115, 2025.

\bibitem{Arr25}
A.~Arroyo-Rabasa.
\newblock Functional and variational aspects of nonlocal operators associated
  with linear {PDE}s.
\newblock {\em Nonlinear Anal.}, 251:Paper No. 113683, 26, 2025.

\bibitem{Ball1977}
J.~M. Ball.
\newblock Convexity conditions and existence theorems in nonlinear elasticity.
\newblock {\em Arch. Rational Mech. Anal.}, 63(4):337--403, 1976/77.

\bibitem{BCFR24}
J.~C. Bellido, J.~Cueto, M.~D. Foss, and P.~Radu.
\newblock Nonlocal {G}reen theorems and {H}elmholtz decompositions for
  truncated fractional gradients.
\newblock {\em Appl. Math. Optim.}, 90(1):Paper No. 16, 49, 2024.

\bibitem{BeCuMC}
J.~C. Bellido, J.~Cueto, and C.~Mora-Corral.
\newblock Fractional {P}iola identity and polyconvexity in fractional spaces.
\newblock {\em Ann. Inst. H. Poincar\'e{} C Anal. Non Lin\'eaire},
  37(4):955--981, 2020.

\bibitem{BCM21}
J.~C. Bellido, J.~Cueto, and C.~Mora-Corral.
\newblock {$\Gamma $}-convergence of polyconvex functionals involving
  {$s$}-fractional gradients to their local counterparts.
\newblock {\em Calc. Var. Partial Differential Equations}, 60(1):Paper No. 7,
  29, 2021.

\bibitem{BeCuMC23}
J.~C. Bellido, J.~Cueto, and C.~Mora-Corral.
\newblock Non-local gradients in bounded domains motivated by continuum
  mechanics: fundamental theorem of calculus and embeddings.
\newblock {\em Adv. Nonlinear Anal.}, 12(1):Paper No. 20220316, 48, 2023.

\bibitem{BeMCSc24}
J.~C. Bellido, C.~Mora-Corral, and H.~Sch\"onberger.
\newblock Nonlocal gradients: {F}undamental theorem of calculus, {P}oincar\'e{}
  inequalities, and embeddings.
\newblock {\em J. Lond. Math. Soc. (2)}, 112(2):Paper No. e70277, 2025.

\bibitem{Comi2}
E.~Bru\`e, M.~Calzi, G.~E. Comi, and G.~Stefani.
\newblock A distributional approach to fractional {S}obolev spaces and
  fractional variation: asymptotics {II}.
\newblock {\em C. R. Math. Acad. Sci. Paris}, 360:589--626, 2022.

\bibitem{CaR23}
P.~M. Campos and J.~F. Rodrigues.
\newblock Unilateral problems for quasilinear operators with fractional {R}iesz
  gradients.
\newblock {\em Preprint, arXiv:2311.18428}, 2023.

\bibitem{CCM25}
M.~Caponi, A.~Carbotti, and A.~Maione.
\newblock {$H$}-compactness for nonlocal linear operators in fractional
  divergence form.
\newblock {\em Preprint, arXiv:2408.10984}, 2025.

\bibitem{Comi1}
G.~E. Comi and G.~Stefani.
\newblock A distributional approach to fractional {S}obolev spaces and
  fractional variation: existence of blow-up.
\newblock {\em J. Funct. Anal.}, 277(10):3373--3435, 2019.

\bibitem{Con07}
J.~B. Conway.
\newblock {\em A course in functional analysis}, volume~96 of {\em Graduate
  Texts in Mathematics}.
\newblock Springer-Verlag, New York, second edition, 1990.

\bibitem{CuKrSc23}
J.~Cueto, C.~Kreisbeck, and H.~Sch\"onberger.
\newblock A variational theory for integral functionals involving
  finite-horizon fractional gradients.
\newblock {\em Fract. Calc. Appl. Anal.}, 26(5):2001--2056, 2023.

\bibitem{CKS25}
J.~Cueto, C.~Kreisbeck, and H.~Sch\"onberger.
\newblock {$\Gamma $}-convergence involving nonlocal gradients with varying
  horizon: {R}ecovery of local and fractional models.
\newblock {\em Nonlinear Anal. Real World Appl.}, 85:Paper No. 104371, 20,
  2025.

\bibitem{Dac08}
B.~Dacorogna.
\newblock {\em Direct methods in the calculus of variations}, volume~78 of {\em
  Applied Mathematical Sciences}.
\newblock Springer, New York, second edition, 2008.

\bibitem{Delia}
M.~D'Elia, M.~Gulian, H.~Olson, and G.~E. Karniadakis.
\newblock Towards a unified theory of fractional and nonlocal vector calculus.
\newblock {\em Fract. Calc. Appl. Anal.}, 24(5):1301--1355, 2021.

\bibitem{DLSTY19}
M.~D'Elia, X.~Li, P.~Seleson, X.~Tian, and Y.~Yu.
\newblock A review of local-to-nonlocal coupling methods in nonlocal diffusion
  and nonlocal mechanics.
\newblock {\em J. Peridyn. Nonlocal Model.}, 4(1):1--50, 2022.

\bibitem{DGLZ}
Q.~Du, M.~Gunzburger, R.~B. Lehoucq, and K.~Zhou.
\newblock A nonlocal vector calculus, nonlocal volume-constrained problems, and
  nonlocal balance laws.
\newblock {\em Math. Models Methods Appl. Sci.}, 23(3):493--540, 2013.

\bibitem{DMT22}
Q.~Du, T.~Mengesha, and X.~Tian.
\newblock Fractional {H}ardy-type and trace theorems for nonlocal function
  spaces with heterogeneous localization.
\newblock {\em Anal. Appl. (Singap.)}, 20(3):579--614, 2022.

\bibitem{DTWY22}
Q.~Du, X.~Tian, C.~Wright, and Y.~Yu.
\newblock Nonlocal trace spaces and extension results for nonlocal calculus.
\newblock {\em J. Funct. Anal.}, 282(12):Paper No. 109453, 63, 2022.

\bibitem{Fef73}
C.~Fefferman.
\newblock {$L\sp{p}$} bounds for pseudo-differential operators.
\newblock {\em Israel J. Math.}, 14:413--417, 1973.

\bibitem{FoL07}
I.~Fonseca and G.~Leoni.
\newblock {\em Modern methods in the calculus of variations: {$L^p$} spaces}.
\newblock Springer Monographs in Mathematics. Springer, New York, 2007.

\bibitem{GHS21}
H.~F. Gon\c{c}alves, D.~D. Haroske, and L.~Skrzypczak.
\newblock Compact embeddings in {B}esov-type and {T}riebel-{L}izorkin-type
  spaces on bounded domains.
\newblock {\em Rev. Mat. Complut.}, 34(3):761--795, 2021.

\bibitem{Gra14b}
L.~Grafakos.
\newblock {\em Modern {F}ourier analysis}, volume 250 of {\em Graduate Texts in
  Mathematics}.
\newblock Springer, New York, third edition, 2014.

\bibitem{HMT24}
Z.~Han, T.~Mengesha, and X.~Tian.
\newblock Compactness results for a {D}irichlet energy of nonlocal gradient
  with applications.
\newblock {\em Numer. Methods Partial Differential Equations}, 40(6):Paper No.
  e23149, 46, 2024.

\bibitem{HaT23}
Z.~Han and X.~Tian.
\newblock Nonlocal half-ball vector operators on bounded domains:
  {P}oincar\'e{} inequality and its applications.
\newblock {\em Math. Models Methods Appl. Sci.}, 33(12):2507--2556, 2023.

\bibitem{Hor07}
L.~H\"{o}rmander.
\newblock {\em The analysis of linear partial differential operators. {III}}.
\newblock Classics in Mathematics. Springer, Berlin, 2007.
\newblock Pseudo-differential operators, Reprint of the 1994 edition.

\bibitem{KrS22}
C.~Kreisbeck and H.~Sch\"{o}nberger.
\newblock Quasiconvexity in the fractional calculus of variations:
  {C}haracterization of lower semicontinuity and relaxation.
\newblock {\em Nonlinear Anal.}, 215:Paper No. 112625, 2022.

\bibitem{KrS24}
C.~Kreisbeck and H.~Sch\"onberger.
\newblock Non-constant functions with zero nonlocal gradient and their role in
  nonlocal {N}eumann-type problems.
\newblock {\em Nonlinear Anal.}, 249:Paper No. 113642, 28, 2024.

\bibitem{Kum81}
H.~Kumano-go.
\newblock {\em Pseudodifferential operators}.
\newblock MIT Press, Cambridge, Mass.-London, 1981.
\newblock Translated from the Japanese by the author, R\'{e}mi Vaillancourt and
  Michihiro Nagase.

\bibitem{MeS}
T.~Mengesha and D.~Spector.
\newblock Localization of nonlocal gradients in various topologies.
\newblock {\em Calc. Var. Partial Differential Equations}, 52(1-2):253--279,
  2015.

\bibitem{Miy87}
A.~Miyachi.
\newblock Estimates for pseudodifferential operators with exotic symbols.
\newblock {\em J. Fac. Sci. Univ. Tokyo Sect. IA Math.}, 34(1):81--110, 1987.

\bibitem{Morrey}
C.~B. Morrey, Jr.
\newblock Quasi-convexity and the lower semicontinuity of multiple integrals.
\newblock {\em Pacific J. Math.}, 2:25--53, 1952.

\bibitem{RaR07}
V.~S. Rabinovich and S.~Roch.
\newblock Exact and numerical inversion of pseudo-differential operators and
  applications to signal processing.
\newblock In {\em Modern trends in pseudo-differential operators}, volume 172
  of {\em Oper. Theory Adv. Appl.}, pages 259--277. Birkh\"{a}user, Basel,
  2007.

\bibitem{Rindler}
F.~Rindler.
\newblock {\em Calculus of variations}.
\newblock Universitext. Springer, Cham, 2018.

\bibitem{Ryc99}
V.~S. Rychkov.
\newblock On restrictions and extensions of the {B}esov and
  {T}riebel-{L}izorkin spaces with respect to {L}ipschitz domains.
\newblock {\em J. London Math. Soc. (2)}, 60(1):237--257, 1999.

\bibitem{ScD24}
J.~M. Scott and Q.~Du.
\newblock Nonlocal problems with local boundary conditions {I}: {F}unction
  spaces and variational principles.
\newblock {\em SIAM J. Math. Anal.}, 56(3):4185--4222, 2024.

\bibitem{ScD25}
J.~M. Scott and Q.~Du.
\newblock Nonlocal problems with local boundary conditions {II}: {G}reen's
  identities and regularity of solutions.
\newblock {\em SIAM J. Math. Anal.}, 57(1):404--451, 2025.

\bibitem{ShY23}
Z.~Shi and L.~Yao.
\newblock New estimates of {R}ychkov's universal extension operator for
  {L}ipschitz domains and some applications.
\newblock {\em Math. Nachr.}, 297(4):1407--1443, 2024.

\bibitem{ShS2015}
T.-T. Shieh and D.~E. Spector.
\newblock On a new class of fractional partial differential equations.
\newblock {\em Adv. Calc. Var.}, 8(4):321--336, 2015.

\bibitem{ShS2018}
T.-T. Shieh and D.~E. Spector.
\newblock On a new class of fractional partial differential equations {II}.
\newblock {\em Adv. Calc. Var.}, 2017.

\bibitem{Silhavy2019}
M.~{\v{S}}ilhav\'y.
\newblock Fractional vector analysis based on invariance requirements (critique
  of coordinate approaches).
\newblock {\em Contin. Mech. Thermodyn.}, 32(1):207--228, 2020.

\bibitem{Sil00}
S.~A. Silling.
\newblock Reformulation of elasticity theory for discontinuities and long-range
  forces.
\newblock {\em J. Mech. Phys. Solids}, 48(1):175--209, 2000.

\bibitem{SEWXA07}
S.~A. Silling, M.~Epton, O.~Weckner, J.~Xu, and E.~Askari.
\newblock Peridynamic states and constitutive modeling.
\newblock {\em J. Elasticity}, 88(2):151--184, 2007.

\bibitem{SLS14}
S.~A. Silling, D.~J. Littlewood, and P.~Seleson.
\newblock Variable horizon in a peridynamic medium.
\newblock {\em J. Mech. Mater. Struct.}, 10(5):591--612, 2015.

\bibitem{Stein}
E.~M. Stein.
\newblock {\em Singular integrals and differentiability properties of
  functions}, volume No. 30 of {\em Princeton Mathematical Series}.
\newblock Princeton University Press, Princeton, NJ, 1970.

\bibitem{TTD19}
Y.~Tao, X.~Tian, and Q.~Du.
\newblock Nonlocal models with heterogeneous localization and their application
  to seamless local-nonlocal coupling.
\newblock {\em Multiscale Model. Simul.}, 17(3):1052--1075, 2019.

\bibitem{Tay11b}
M.~E. Taylor.
\newblock {\em Partial differential equations {III}. {N}onlinear equations},
  volume 117 of {\em Applied Mathematical Sciences}.
\newblock Springer, New York, second edition, 2011.

\bibitem{TiD17}
X.~Tian and Q.~Du.
\newblock Trace theorems for some nonlocal function spaces with heterogeneous
  localization.
\newblock {\em SIAM J. Math. Anal.}, 49(2):1621--1644, 2017.

\bibitem{Trie1}
H.~Triebel.
\newblock {\em Theory of function spaces}, volume~78 of {\em Monographs in
  Mathematics}.
\newblock Birkh\"{a}user Verlag, Basel, 1983.

\end{thebibliography}
\end{document}